\documentclass[a4paper,11pt]{article}
\usepackage{graphicx}
\usepackage{amsmath,amssymb,amsthm}
\usepackage{mathtools}
\usepackage{hyperref}
\usepackage{breakurl}
\usepackage{mhequ}
\usepackage{mathrsfs}
\usepackage{ifthen}
\usepackage{tikz}
\usepackage{microtype}
\usepackage{wasysym}
\usepackage{centernot}
\usepackage{enumitem}
\usepackage{shuffle}
\usepackage{natbib}
\usepackage{pgfplots}
\usepackage{tikz-3dplot}
\pgfplotsset{compat=newest}

\usetikzlibrary{math}

\usepackage{mathtools}

\usetikzlibrary{shapes,
                chains,
                scopes}

\usepackage{geometry}
\usetikzlibrary{shapes.misc}
\usetikzlibrary{shapes.symbols}
\usetikzlibrary{shapes.geometric}
\usetikzlibrary{decorations}
\usetikzlibrary{decorations.markings}
\usetikzlibrary{patterns} 
\usetikzlibrary{snakes} 

\usepackage{caption}
\usepackage{subcaption}
\usetikzlibrary{plotmarks}
\usetikzlibrary{decorations.pathreplacing}
\usetikzlibrary{decorations.pathmorphing}
\usetikzlibrary{calc}
\usetikzlibrary{external}
\usetikzlibrary{arrows.meta}

\usepgfplotslibrary{fillbetween}

\theoremstyle{plain}
\newtheorem{theorem}{Theorem}[section]
\newtheorem{lemma}[theorem]{Lemma}

\theoremstyle{definition}
\newtheorem{definition}[theorem]{Definition}

\newtheorem{example}[theorem]{Example}

\def\bone{\mathbf{1}}

\newcommand{\bbX}{\mathbb{X}}
\def\R{\mathbb{R}}
\def\E{\mathbb{E}}

\def\Z{\mathbb{Z}}

\newcommand{\mfT}{\mathfrak{T}}
\newcommand{\mfn}{\mathfrak{n}}

\newcommand{\mfg}{\mathfrak{g}}

\newcommand{\mbY}{\mathbf{Y}}
\newcommand{\mbX}{\mathbf{X}}
\newcommand{\mbx}{\mathbf{x}}

\newcommand{\mrd}{\mathop{}\!\mathrm{d}}

\newcommand{\id}{\mathrm{id}}

\newcommand\scal[2][ ]{\ifthenelse{\equal{#1}{ }}{\langle#2\rangle}{}
        \ifthenelse{\equal{#1}{b}}{\bigl\langle#2\bigr\rangle}{}
        \ifthenelse{\equal{#1}{B}}{\Bigl\langle#2\Bigr\rangle}{}
        \ifthenelse{\equal{#1}{bb}}{\biggl\langle#2\biggr\rangle}{}
        \ifthenelse{\equal{#1}{BB}}{\Biggl\langle#2\Biggr\rangle}{}}
\newcommand{\eps}{\varepsilon}

\newcommand{\CK}{\mathrm{CK}}
\newcommand{\GL}{\mathrm{GL}}
\newcommand{\Lip}{\mathrm{Lip}}

\newcommand{\floor}[1]{\lfloor #1 \rfloor}
\newcommand{\var}[1]{#1\textnormal{-var}}
\newcommand{\Hol}[1]{#1\textnormal{-H{\"o}l}}
\newcommand{\grp}{G\Omega}
\newcommand{\wgrp}{W\Omega}
\newcommand{\Span}{\mathrm{span}}

\def\dash{\leavevmode\unskip\kern0.18em--\penalty\exhyphenpenalty\kern0.18em}
\def\slash{\leavevmode\unskip\kern0.15em/\penalty\exhyphenpenalty\kern0.15em}

\let\OLDthebibliography\thebibliography
\renewcommand\thebibliography[1]{
  \OLDthebibliography{#1}
  \setlength{\parskip}{0pt}
  \setlength{\itemsep}{2mm}
}

\colorlet{symbols}{blue!90!black}
\colorlet{testcolor}{green!60!black}
\colorlet{connection}{red!30!black}

\tikzset{
root/.style={circle,fill=black!50,inner sep=0pt, minimum size=3mm},
        dot/.style={circle,fill=black,inner sep=0pt, minimum size=1.5mm},
        dotred/.style={circle,fill=black!50,inner sep=0pt, minimum size=2mm},
        var/.style={circle,fill=black!10,draw=black,inner sep=0pt, minimum size=3mm},
        kernel/.style={semithick,shorten >=2pt,shorten <=2pt},
        kernel1/.style={thick},
        kernels/.style={snake=zigzag,shorten >=2pt,shorten <=2pt,segment amplitude=1pt,segment length=4pt,line before snake=2pt,line after snake=5pt,},
		kernels1/.style={snake=zigzag,segment amplitude=0.5pt,segment length=2pt},
		rho1/.style={dotted,semithick},
        rho/.style={densely dashed,semithick,shorten >=2pt,shorten <=2pt},
           testfcn/.style={dotted,semithick,shorten >=2pt,shorten <=2pt},
        renorm/.style={shape=circle,fill=white,inner sep=1pt},
        labl/.style={shape=rectangle,fill=white,inner sep=1pt},
        xic/.style={very thin,circle,fill=symbols,draw=black,inner sep=0pt,minimum size=1.2mm},
        xi/.style={very thin,circle,fill=blue!10,draw=black,inner sep=0pt,minimum size=1.2mm},
        xix/.style={crosscircle,fill=blue!10,draw=black,inner sep=0pt,minimum size=1.2mm},
	xib/.style={very thin,circle,fill=blue!10,draw=black,inner sep=0pt,minimum size=1.6mm},
	xie/.style={very thin,circle,fill=green!50!black,draw=black,inner sep=0pt,minimum size=1mm},
	xid/.style={very thin,circle,fill=symbols,draw=black,inner sep=0pt,minimum size=1.6mm},
	xibx/.style={crosscircle,fill=blue!10,draw=black,inner sep=0pt,minimum size=1.6mm},
	edgetype/.style={very thin,circle,draw=black,inner sep=0pt,minimum size=5mm},
	nodetype/.style={very thick,circle,draw=black,inner sep=0pt,minimum size=5mm},
	kernels2/.style={very thick,draw=connection,segment length=12pt},
clean/.style={thin,circle,fill=black,inner sep=0pt,minimum size=1mm},	not/.style={thin,circle,fill=symbols,draw=connection,fill=connection,inner sep=0pt,minimum size=0.5mm},
	>=stealth,
        }

\pgfplotsset{soldot/.style={color=blue,only marks,mark=*}} \pgfplotsset{holdot/.style={color=blue,fill=white,only marks,mark=*}}
\pgfplotsset{soldotred/.style={color=red,only marks,mark=*}} \pgfplotsset{holdotred/.style={color=red,fill=white,only marks,mark=*}}

\providecommand{\keywords}[1]
{
  \small	
  \textbf{\textit{Keywords---}} #1
}

\title{Rough path theory}

\author{
Ilya Chevyrev\thanks{School of Mathematics, The University of Edinburgh, Edinburgh EH9 3FD, United Kingdom.
\href{mailto:ichevyrev@gmail.com}{\tt ichevyrev@gmail.com}}
}

\date{}

\begin{document}

\maketitle

\begin{abstract}
The theory of rough paths arose from a desire to establish continuity properties of ordinary differential equations involving terms of low regularity. While essentially an analytic theory, its main motivation and applications are in stochastic analysis, where it has given a new perspective on It\^o calculus and a meaning to stochastic differential equations driven by irregular paths outside the setting of semi-martingales. In this survey, we present some of the main ideas that enter rough path theory. We discuss complementary notions of solutions for rough differential equations and the related notion of path signature, and give several applications and generalisations of the theory.
\end{abstract}

\keywords{Rough paths, stochastic analysis, sewing lemma, iterated integrals, path signature, stochastic partial differential equations, regularity structures}

\tableofcontents

\section{Introduction}\label{sec:intro}

Classical It\^o calculus is able to make sense of stochastic differential equations (SDEs) of the form
\[
\mrd Y_t = f(Y_t)\mrd X_t\;,\quad Y_0\in\R^n\;,
\]
where $X\colon [0,T]\to \R^d$ is a semi-martingale, such as  Brownian motion (Wiener process).
A fundamental structure that It\^o calculus exploits is the martingale property, and the solution $Y$ is inherently a probabilistic object.
For example, the It\^o integral is defined as a random variable with reference to the entire law of $X$ rather than to a single realisation.
In particular, while the solution map $X\mapsto Y$ is measurable, it is not continuous,
and thus, given a single trajectory of $X$,
there is no \textit{canonical} way to construct from it a solution $Y$.
This is also evident from the fact that there are multiple ways to interpret the product $f(Y)\mrd X$, e.g. in the sense of It\^o or Stratonovich,
both of which can be `correct' depending on the context and which in general give different answers.

The theory of rough paths was introduced by Lyons in the 90's as a way to understand and overcome this ambiguity and lack of continuity in SDEs.
Lyons observed that, while there is no way to extend Riemann--Stieltjes--Young integration in a continuous manner to Brownian motion using linear functional analysis,
there is essentially only one non-linear object that needs to be chosen, the 2nd iterated integral,
which then removes all ambiguity in the notion of solution to SDEs and renders the solution map continuous.
This leads to the notion of a rough path $\mbX$, which can be seen as an underlying path $X$ (e.g. a trajectory of Brownian motion)
together with a finite number of higher order objects $\mbX^2,\mbX^3,\ldots,\mbX^N$ that are \textit{a priori} not canonically defined just from $X$
(for Brownian motion, only the 2nd iterated integral $\mbX^2 = \int X\mrd X$ is needed).
These higher order objects cannot be chosen freely and are required to respect important algebraic and analytic constraints.
 
The space of rough paths is non-linear and comes equipped with a $p$-variation or H\"older-type metric
that is stronger than the analogous metric at the level of the path $X$.
Crucially, one can extend to the space of rough paths the solution map $\mbX\mapsto Y$,
which is now continuous under the stronger rough path metric.
By `extend', we mean that there is a canonical injection of smooth paths into rough paths and that the classical solution map for smooth paths factors through this injection and the rough path solution map.

SDEs can thus be solved, for example, in the Stratonovich sense
by choosing once and for all the 2nd iterated integral $\int X\mrd X$ of Brownian motion $X$ in the Stratonovich sense.
This gives rise to a random rough path $\mbX$ from which the Stratonovich solution $Y$ to any SDE is derived in a \textit{pathwise} and continuous manner.
A different choice for the 2nd iterated integral (e.g. in the It\^o sense), will give rise a different random rough path $\tilde \mbX$
and thus different solutions $\tilde Y$ to SDEs but which are still continuous functions of $\tilde \mbX$.
This continuity further leads to approximation results, e.g. if $X^n$ are smooth paths such that $X^n\to X$ and $\int X^n\mrd X^n \to \int X \mrd X$ in the rough path metric,
then the associated solutions $Y^n$ also converge to $Y$.

Rough path theory therefore factors the solution map for SDEs into two steps: (1) the (typically probabilistic) construction of a rough path $\mbX$ from a stochastic process $X$,
and (2) a deterministic and continuous solution map.
Importantly, this separation of probability and analysis has allowed rough path theory to move beyond the setting of semi-martingales.
For example, rough paths can be built in a canonical way from fractional Brownian motion with Hurst parameter $H\in(\frac14,\frac12)$, which is not a semi-martingale
and falls outside the scope of It\^o calculus.
In fact, rough paths can be studied from an entirely analytic perspective without the need for an underlying stochastic processes.

In this survey, we aim to expose the key ideas that enter rough path theory.
While several monographs exist on the subject, see in particular \citep{LyonsQian02, FV10, FrizHairer20}, we aim for a complementary exposition and discuss several aspects that are not covered in detail by any one of these works.
This includes Gubinelli's approach to rough path theory via controlled rough paths,
Davie's notion of solution based on Euler estimates,
extensions of rough path and signatures,
generalisations to branched (non-geometric) and discontinuous rough paths,
and a brief discussion of how rough paths extend to regularity structures.
Important topics that we do not discuss include interactions with Malliavin calculus,
the stochastic sewing lemma,
Doob--Meyer decomposition and true roughness,
and rough partial differential equations.

\section{Sewing and Young integration}

One of the core ideas behind rough path theory is the following sewing lemma.
Let $\Delta_T = \{(s,t)\,:\, 0\leq s\leq t\leq T\}$.
A \textit{control function} is a continuous map $\omega\colon \Delta_T \to [0,\infty)$ which is super-additive, i.e. $\omega(s,t)+\omega(t,u)\leq\omega(s,u)$ for all $0\leq s\leq t\leq u\leq T$.

\begin{lemma}[Sewing]
Consider $\beta>1$, a real Banach space $V$, and a control function $\omega$. Let $\Xi\colon \Delta_T\to V$ be a continuous map such that, for all $0\leq s\leq u\leq t\leq T$,
\[
\|\delta\Xi_{s,u,t}\| \leq \omega(s,t)^\beta\;,
\]
where $\delta\Xi_{s,u,t} = \Xi_{s,t}-\Xi_{s,u}-\Xi_{u,t}$.
Then there exists a unique function $\xi\colon [0,T]\to V$ such that $\xi_0=0$ and
\[
|\xi_{s,t}-\Xi_{s,t}| \leq C \omega(s,t)^\beta\;,
\]
where we denote $\xi_{s,t}=\xi_t-\xi_s$
and $C$ is a constant depending only on $\beta$.
Furthermore, the increments of $\xi$ are given by the limit of Riemann sums $\xi_{s,t} = \lim_{|P|\to 0}\sum_{[u,v]\in P} \Xi_{u,v}$.
\end{lemma}

Above, $P$ is a dissection of $[s,t]$, by which we mean a finite partition of $[s,t]$ into disjoint intervals of strictly positive length with mesh size $|P|:=\max_{[u,v]\in P}|u-v|$, and where we write $[u,v]\in P$ to mean that $u<v$ are the endpoints of an interval in $P$.
This result (in the case $\omega(s,t)=|t-s|$)
is due to \cite{Gubinelli04} and \cite{Feyel_de_La_Pradelle_06};
the above version is from \citep{DGHT19}.

The picture to keep in mind is that $\delta\Xi$ measures how far $\Xi$ is from being the increments of a path -- if $\delta\Xi=0$, then evidently the unique $\xi$ is $\xi_{t} = \Xi_{0,t}$ and $\Xi_{s,t}=\xi_{s,t}$.
The point of the lemma is that it requires $\delta\Xi$ to only be small in a suitable sense and this suffices to recover a path $\xi$ whose increments $\Xi$ approximates.
A proof of an analogous result dates back to L. C. Young in 1936.
There exists an extension of the sewing lemma to $\beta\leq 1$ due to \cite{BZ22}, wherein one loses uniqueness of $\xi$ as well as its construction 
as a limit of Riemann sums.

A typical application of the sewing lemma is to show that, for $X\in C^{\var p}([0,T],V)$ and $Y\in C^{\var q}([0,T],L(V,W))$ where $\frac1p+\frac1q>1$ with $p,q\in[1,\infty)$, the so-called Young integral $\int_0^T Y_u\mrd X_u$ is well-defined as a limit of Riemann sums.
Here $V,W$ are Banach spaces, $L(V,W)$ is the space of bounded linear operators from $V$ to $W$
equipped with the operator norm,
and $C^{\var p}([a,b],V)$ is the vector space of continuous functions $Z\colon [a,b]\to V$ equipped with $p$-variation semi-norm
\begin{equation}\label{eq:p_var_linear}
\|Z\|_{\var p;[a,b]}:=\omega_{Z,p}(a,b)^{1/p}:=\sup_{P} \Big(\sum_{[s,t]\in P} \|Z_{s,t}\|^p\Big)^{1/p}\;,
\end{equation}
where the supremum is over all dissections $P$ of $[a,b]$.

To construct the integral $\int_0^T Y_u\mrd X_u$, let us define $\Xi_{s,t}=Y_sX_{s,t}$ and observe that
\[
\delta\Xi_{s,u,t} = Y_{s,u}X_{t,u} \leq \|Y\|_{\var q;[s,t]}\|X\|_{\var p;[s,t]}= \omega(s,t)^{1/p+1/q}\;,
\]
where $\omega(s,t)=(\omega_{Y,q}^p\omega_{X,p}^q)^{1/(p+q)}(s,t)$ is a control function since both $\omega_{Y,q}$ and $\omega_{X,p}$ are.
Hence, by the sewing lemma, there exists a unique additive function
\begin{equation}\label{eq:Y_X_int}
(s,t)\mapsto \int_s^t Y_u\mrd X_u \in W\;,
\end{equation}
given by a limit of Riemann sums,
for which
\[
\Big|\int_s^t Y_u\mrd X_u-Y_sX_{s,t}\Big|
\leq C\|X\|_{\var p;[s,t]}\|Y\|_{\var q;[s,t]}\;.
\]
The function $\int_0^\bullet Y_u\mrd X_u$ is furthermore in $C^{\var p}([0,T],V)$
and is jointly continuous in $(X,Y)$
once we equip the state space of $Y$ with the norm $\|Y_0\|+\|Y\|_{\var q;[0,T]}$.

If $\|X\|_{\var 1;[0,T]}<\infty$, then $X$ is said to be a bounded variation.
In this case the integral \eqref{eq:Y_X_int} is the classical Riemann--Stieltjes integral (applicable also to $q=\infty$, i.e. $Y\in C([0,T],L(V,W))$ equipped with the uniform norm).

By using the control function $\omega(s,t)=|t-s|$ in the sewing lemma, one recovers the same result specialised to the H\"older scale, namely that $(X,Y)\mapsto \int_0^\bullet Y_u\mrd X_u$ is well defined and jointly continuous on the space $C^\alpha\times C^\beta \to C^\alpha$ if $\alpha+\beta>1$,
where $C^\alpha$ is equipped with the H\"older semi-norm
\[
\|Z\|_{\Hol\alpha} = \sup_{s\neq t} |t-s|^{-\alpha}\|Z_{s,t}\|
\] 
and $C^\beta$ with the norm $\|Z\|_{C^\beta} = \|Z_0\| + \|Z\|_{\Hol\beta}$.

Coming back to differential equations, by contraction mapping arguments,
one can show that the differential equation
\begin{equation}\label{eq:ODE}
\mrd Y_t = f(Y_t)\mrd X_t\;,\quad Y_0 \in W
\end{equation}
admits a unique solution,
by which we mean $Y_t = Y_0 + \int_0^t f(Y_s)\mrd X_s$ for all $t\in[0,T]$,
whenever $X\in C^{\var p}([0,T],V)$ for $p\in [1,2)$
and $f\colon W\to L(V,W)$ is sufficiently smooth with bounded derivatives ($f\in C^2$ suffices).
Furthermore the solution map $X\mapsto Y$ is locally Lipschitz.
By using the control function $\omega(s,t)=|t-s|$ as earlier, one recovers the same result in the H\"older scale $C^\alpha$ for $\alpha\in (\frac12,1]$.

As a first application to stochastic analysis, consider a $d$-dimensional fractional Brownian motion
\[
X=(X^1,\ldots, X^d)\colon[0,T]\to\R^d
\]
with Hurst parameter $H\in (0,1)$, i.e.
$X$ is a centred Gaussian process with $X^i,X^j$
independent for $i\neq j$
and
\[
\E X^i_sX^i_t = \frac12 (t^{2H}+s^{2H}-|t-s|^{2H})\;.
\]
If $H=\frac12$, then $X$ is a classical Brownian motion.
By an argument due to Kolmogorov, one can show that almost every trajectory of $X$ belongs to the H\"older space $C^{\alpha}$ for any $\alpha<H$.
Therefore, in the regime $H>\frac12$,
the differential equation \eqref{eq:ODE}, which is now a \textit{stochastic differential equation} (SDE),
admits a unique solution.
Furthermore, the continuity of the solution map $X\mapsto Y$ readily implies a number of results about the solution; for example, it implies the almost sure convergence $Y^\eps \to Y$,
where $Y^\eps$ solves \eqref{eq:ODE} with $X$ replaced by a smooth approximation $X^\eps$  such that $X^\eps \to X$ in $C^\alpha$ almost surely, $\alpha>\frac12$ (e.g. mollification or piecewise linear interpolation),
which is a multidimensional and fractional version of the classical Wong--Zakai theorem \citep[see][Chap.~VI, Thm.~7.2]{Ikeda_Watanabe_89_SDEs}.

One should contrast this solution theory to It\^o SDEs where one relies on a martingale structure of $X$, thus applicable only to $H=\frac12$,
and which lacks pathwise continuity of the map $X\mapsto Y$.
On the other hand, Young integration is unable to handle $H\leq \frac12$, in particular $H=\frac12$, which is arguably the most important
case in stochastic analysis.

To see the difficulty of moving from $\alpha>\frac12$ to $\alpha\leq \frac12$ in \eqref{eq:ODE},
it is instructive to note that the iterated integral $\int_0^T X^1_t\mrd X^2_t$ for a $2$-dimensional path $X\colon [0,T]\to \R^2$ is a basic example of a solution to \eqref{eq:ODE} (with $f$ linear).
For $\alpha>\frac12$, the map $X\mapsto \int_0^T X^1_t\mrd X^2_t$, defined for smooth $X$, extends continuously to $C^\alpha$, while for $\alpha\leq \frac12$ this map does not have a continuous extension to $C^\alpha$ and in fact
does not even have a closable graph, i.e.
there exists a sequence of smooth paths $X(n)\colon [0,T]\to \R^2$ such that $\|X(n)\|_{\Hol\alpha}\to 0$ while $\int_0^T X^1_t(n)\mrd X^2_t(n) = 1$ for all $n$.
The next example shows this for $\alpha<\frac12$.

\begin{example}
Consider $X\in C^\infty([0,2\pi],\R^2)$ given by
$X_t = n^{-1}(\cos(n^2t),\sin(n^2t))$.
Then $\|X(n)\|_{\Hol\alpha} = \pi^\alpha n^{2\alpha-1} \to 0$ as $n\to\infty$ for any $\alpha<\frac12$, but
\[
\int_0^{2\pi} X^1_t(n)\mrd X^2_t(n) =  \int_0^{2\pi}\cos^2(n^2t)\mrd t=\pi\;.
\]
This is can be understood pictorially: let $A^{i,j}_t = \int_0^t X^i_s\mrd X^j_s$
and observe that the so-called L\'evy area $\frac12(A^{1,2}_{t}-A^{2,1}_{t})$ is the signed area traced out by $(X^1,X^2)$ and the chord connecting its endpoints, see Figure \ref{fig:LevyArea}.
Then for $X(n)$, the L\'evy area up to time $t$ is $t/2$, see Figure \ref{fig:spiral}.
\end{example}

\begin{figure}
\begin{subfigure}[h]{0.48\textwidth}
\centering
\begin{tikzpicture}[thick,scale=0.8,
declare function={ f(\x)=0.9*(\x-1.6)^3 - (\x-1.6) + 1.1; 
}]

\begin{axis}[
    width=8cm, 
    height=8cm, 
    ticks=none]
\addplot[{Latex[length=3mm]}-,
    name path=A, 
    mark=none,
    line width=2pt,
    domain=0.2:2.8,
    samples=201
    ]
    {f(x)};
    
\addplot[
    name path=B,
    mark=none,
    dashed,
    ]
    coordinates {(0.2,{f(0.2)}) (2.8,{f(2.8)})};

\addplot[] fill between [of=A and B,split,every segment no 0/.style={color=gray!40},every segment no 1/.style={color=red!0},];

\node at (axis cs: 1, 1.0) {\huge$\mathbf{+}$};
\node at (axis cs: 2.25, 0.9) {\huge$\mathbf{-}$};
\end{axis}
\end{tikzpicture}
\caption{L\'evy area of a two-dimensional path (grey area minus white area).}\label{fig:LevyArea}
\end{subfigure}
\hfill
\begin{subfigure}[h]{0.48\textwidth}
\centering
\begin{tikzpicture}[scale=0.8]
\begin{axis}[
    width=8cm, 
    height=8cm, 
    view={-20}{-20},
    axis line style = ultra thick,
    axis lines=middle,
    zmax=80,
    xmax=2,
    ymax=2,
    xtick=\empty,
    ytick=\empty,
    ztick=\empty,
    clip=false,
    x label style={at={(axis cs:2,0.051)},anchor=north},
    xlabel={$y$},
    y label style={at={(axis cs:0.05,2)},anchor=north},
    ylabel={$x$},
    z label style={at={(axis cs:0.075,0,80)},anchor=north},
]
\draw[ultra thick,red] (0,0) circle[radius=1];
\addplot3+[domain=0:11*pi,samples=500,samples y=0,black,no marks,ultra thick] 
({sin(deg(x))}, 
{cos(deg(x))}, 
{6*x/(pi)})
node [circle,scale=0.5,fill,pos=0]{};
\end{axis}
\end{tikzpicture}
\caption{The path $X(n)$ with its L\'evy area indicated in the $z$-axis (circle in $x$-$y$ plane has radius $\frac1n$).}\label{fig:spiral}
\end{subfigure}
\caption{}
\end{figure}

This lack of continuity of iterated integrals persists at a probabilistic level when considering Brownian motion, and is independent of the choice of norm used to measure regularity (e.g. $\alpha$-H\"older or $p$-variation).
This is made precise by the following result of \cite{Lyons91}: there exists no Banach space $V\subset C([0,1],\R)$ such that (a) all smooth paths are in $V$, (b) Fourier approximations of Brownian motion converge in $V$, and (c) the bilinear map $(X,Y)\mapsto \int_0^1 Y_t\mrd X_t$, defined for smooth $X,Y$, extends to a continuous bilinear map on $V$.
This shows, in essence, that there is no linear integration theory able to extend Riemann integration to the case of Brownian motion in a continuous way.

\section{Iterated integrals}

Iterated integrals above provided an obstruction to the continuity of $X\mapsto Y$ for \eqref{eq:ODE} in the regime $X\in C^{\var p}$ with $p\geq 2$.
It turns out that these iterated integrals are the \textit{only} obstruction in the following sense:
once we make a sensible definition for these iterated integrals, and enhance the path $X$ to a richer object $\mbX$ that contains this new data,
one can recover continuity of a richer map $\mbX\mapsto Y$ that extends the classical solution map (see the Universal Limit Theorem \ref{thm:ULT} below).
This leads to the following definition.

\begin{definition}[Rough path]\label{def:RP}
For $p\in[1,\infty)$ and a Banach space $V$, a $p$-rough path
is a map
\[
\mbX=(\mbX^0,\ldots, \mbX^{N}) \colon \Delta_T \to T^{N}(V)\;,
\]
where $N=\floor{p}$ is the floor of $p$ and $T^{N}(V)$ is the truncated tensor algebra
\[
T^{N}(V)= \bigoplus_{k=0}^{N}V^{\otimes k}\;,
\]
which satisfies
\begin{enumerate}[label=(\alph*)]
\item\label{pt:Chen} Chen's identity: $\mbX_{s,t} \mbX_{t,u}=\mbX_{s,u}$ for all $s\leq t \leq u$,
\item\label{pt:pvar} finite $p$-variation: there exists a control function $\omega$ such that, for all $(s,t)\in\Delta_T$,
\begin{equation}\label{eq:p-var-omega}
\max_{1\leq k\leq N}\|\mbX^k_{s,t}\|_{V^{\otimes k}}^{p/k} \leq \omega(s,t)\;,
\end{equation}
\item\label{pt:invert} $\mbX^0_{s,t}=1$ for all $(s,t)\in\Delta_T$.
\end{enumerate}
We call $N$ the \textit{level} of $\mbX$.
If one can take $\omega(s,t)=|t-s|$ in \ref{pt:pvar},
then $\mbX$ is called a \textit{$\frac1p$-H\"older rough path}.
\end{definition}

Above, $V^{\otimes k}$ is the completion of the $k$-fold tensor product of $V$ with itself under the projective tensor norm $\|\cdot\|_{V^{\otimes k}}$ (other choices of tensor norms are possible).
By convention $V^{\otimes 0}=\R$.
The product $\mbX_{s,t} \mbX_{t,u}$ is the tensor (i.e. concatenation) product 
\begin{align}\label{eq:tensor_prod}
(x^0,\ldots , x^{N}) (y^0,\ldots , y^{ N})
= (x^0y^0, x^1y^0+x^0y^1,\ldots,
\sum_{i=0}^N x^i y^{N-i})\;.
\end{align}
where $x^iy^j\in V^{\otimes (i+j)}$ is shorthand for $x^i\otimes y^j$.
Below, given $p\in [1,\infty)$, we use the notation $N=\floor p$ unless otherwise stated.

Whenever \ref{pt:pvar} holds, the smallest control function $\omega$ that satisfies \eqref{eq:p-var-omega} is
\begin{equation}\label{eq:p-var}
\omega_{\mbX,p}(s,t) = \sup_{P}  \sum_{[u,v]\in P} \max_{1\leq k\leq N} \|\mbX^k_{u,v}\|_{V^{\otimes k}}^{p/k}\;,
\end{equation}
where the supremum is over all dissections $P$ of $[s,t]$ as in \eqref{eq:p_var_linear}.
The quantity
$\|\mbX\|_{\var p;[s,t]}:=\omega_{\mbX,p}(s,t)^{1/p}$ is called the \textit{$p$-variation} of $\mbX$ over $[s,t]$
and is a generalisation of the semi-norm $\|\cdot\|_{\var p;[s,t]}$ from \eqref{eq:p_var_linear}.

If $X_t-X_s=\mbX^1_{s,t}$ for a path $X\colon [0,T]\to V$, we say that $\mbX$ is a rough path above $X$.
One should think of $\mbX^k_{s,t}$ as a placeholder for the $k$-th iterated integral of $X$.
Indeed,
when $X$ is of bounded variation (or, more generally, finite $q$-variation for $q\in [1,2)$),
then one can define
\begin{equation}\label{eq:II}
\mbX^k_{s,t}=\int_{s}^t\int_s^{t_k}\ldots \int_s^{t_{2}} \mrd X_{t_1} \otimes \ldots \otimes \mrd X_{t_k}
\end{equation}
by means of Riemann--Stieltjes integration (or Young integration), which yields a $p$-rough path for any $p\geq 1$ (or $p\geq q$).
This is called the canonical $p$-rough path lift of $X$.
If $p\geq 2$, this is, however, not the only $p$-rough path above $X$,
e.g. for $p=2$, one can define $\tilde \mbX^2_{s,t} = \mbX^2_{s,t} + Z_t-Z_s$ for any bounded variation $Z\colon[0,T]\to V^{\otimes 2}$
and observe that $(1,\mbX^1,\tilde \mbX^2)$ is again a $p$-rough path above $X$.
Moreover, if $X$ lacks sufficient regularity for the iterated integrals in \eqref{eq:II} to be canonically defined,
then $\mbX^k_{s,t}$ is instead taken as a definition for the right-hand side of \eqref{eq:II} (rather than the other way around).

Condition \ref{pt:Chen} is purely algebraic and is satisfied for $\mbX$ defined by \eqref{eq:II} whenever $X_t-X_s=\mbX^1_{s,t}$ is of bounded variation,
as first observed by Chen.
Condition \ref{pt:pvar}, finiteness of $p$-variation, is analytic and, in light of \eqref{eq:p-var}-\eqref{eq:II},
is the natural extension of $p$-variation of $X$.
Condition \ref{pt:invert} ensures invertibility of $\mbX_{s,t}$ in $T^N(V)$ under the tensor product and, together with \ref{pt:Chen}, implies $\mbX_{t,t}=1$ for all $t\in[0,T]$.
One can in particular view a $p$-rough path equivalently as a function $\mbX\colon[0,T]\to T^N(V)$ whose `increments' $\mbX_{s,t} := \mbX_s^{-1}\mbX_t$ satisfy the conditions in Definition \ref{def:RP}.

\begin{definition}[Geometric rough path]
A $p$-rough path is called \textit{geometric}
if it is the limit of a sequence of smooth paths (lifted to $p$-rough paths via \eqref{eq:II})
under the $p$-variation metric
\begin{equation}\label{eq:d_p-var}
\rho_{\var p}(\mbX,\mbY) := \sup_{P}  \max_{1\leq k\leq N}
\Big(\sum_{[u,v]\in P}  \|\mbX^k_{u,v}-\mbY^k_{u,v}\|_{V^{\otimes k}}^{p/k}\Big)^{k/p}\;,
\end{equation}
where the supremum is over all dissections $P$ of $[0,T]$.
We denote by $\grp_p(V)$ the metric space of geometric $p$-rough paths.
\end{definition}

Considering for now the finite-dimensional case $V=\R^d$,
geometric $p$-rough paths take values in the step-$N$ free nilpotent Lie group $G^N(\R^d)\subset T^N(\R^d)$.
We give two descriptions of the group $G^N(\R^d)$, the first algebraic, the second geometric.

Let $T(\R^d) = \bigoplus_{k=0}^\infty (\R^d)^{\otimes k}$
denote the tensor algebra over $\R^d$.
A basis for $T(\R^d)$ is the set of pure tensors $\{e_{i_1} \cdots  e_{i_k}\}_{i_1,\ldots, i_k}$ for $k \geq 0$ and $1\leq i_1,\ldots, i_k\leq d$ and $\{e_i\}_{i=1}^d$ is the canonical basis of $\R^d$. We call $k$ the length of the tensor $e_{i_1} \cdots  e_{i_k}\in (\R^d)^{\otimes k}$.

We define on $T(\R^d)$ the shuffle product $\shuffle$ by the recursion
\begin{equation}\label{eq:shuffles_def}
(v e_i)\shuffle (we_j) = (v\shuffle (we_j)) e_i + (ve_i\shuffle w)e_j
\end{equation}
for all $v,w\in T(\R^d)$ and $1\leq i,j\leq d$
together with $1\shuffle v = 1\shuffle v = v$
and extended by linearity.
Concretely, $e_{i_1}\cdots e_{i_k}\shuffle e_{j_1}\cdots e_{j_n} = e_{h_1}\cdots e_{h_{k+n}}$ where the multi-index $(h_1,\ldots, h_{k+n})$ is obtained by `shuffling' $(i_1,\ldots, i_{k})$  and $(j_1,\ldots, j_{n})$  while keeping the respective order of terms in the two multi-indexes. For example,
\[
e_{1} e_{2} \shuffle e_{3} = e_{1}e_2e_3 + e_{1}e_3e_2 + e_3e_1e_2\;.
\]
Note that $\shuffle$ turns $T(\R^d)$ in a commutative algebra (in fact, $T(\R^d)$ becomes a Hopf algebra when further equipped with the deconcatenation coproduct).

Next, the space of tensor series $T((\R^d)) := \prod_{k=0}^\infty (\R^d)^{\otimes k}$ is canonically isomorphic to the dual of $T(\R^d)$ (once $\R^d$ is identified with its own dual by taking
 $\{e_i\}_{i=1}^d$ an orthonormal basis of $\R^d$),
and we denote this duality by $\scal{x,y}$ for $x\in T((\R^d))$ and $y\in T(\R^d)$.
We let
\begin{equation}\label{eq:G_def}
G(\R^d) = \{g \in T((\R^d))\,:\, \scal{g,1}=1\,,\, \scal{g,u\shuffle v}=\scal{g,u}\scal{g,v} \; \forall u,v\in T(\R^d) \}
\end{equation}
denote the set of characters on $T(\R^d)$ for the shuffle product.
Then $G(\R^d)$ is a group under the tensor product on $T((\R^d))$.
Let $\pi_N \colon T((\R^d)) \to T^N(\R^d)$ denote the canonical projection and define
\begin{equation}\label{eq:GN}
G^N(\R^d) = \pi_N G(\R^d)\;,
\end{equation}
which is a group in $T^N(\R^d)$ under tensor multiplication.
Equivalently, $G^N(\R^d)$ consists of all $g\in T^N(\R^d)$ such that $\scal{g,u\shuffle v}=\scal{g,u}\scal{g,v}$ for all pure tensors $u,v\in T^N(\R^d)$ of combined length at most $N$.

Another characterisation of $G^N(\R^d)$ is $G^N(\R^d)=\exp(\mfg^N(\R^d))$, where $\mfg^N(\R^d)$ is the free $N$-step nilpotent Lie algebra generated by $\R^d$,
i.e. $\mfg^N(\R^d)$ is the smallest Lie subalgebra of $T^N(\R^d)$ under the tensor product \eqref{eq:tensor_prod} that contains $\R^d$,
and $\exp(X)=\sum_{k=0}^\infty \frac{X^k}{k!}$ is the exponential under the tensor product.
Therefore, $G^N(\R^d)$ is a Carnot group
and, for smooth $X\colon [0,T]\to\R^d$, the geometric $p$-rough path $[0,T]\ni t\mapsto \mbX_{0,t}\in G^N(\R^d)$
defined by \eqref{eq:II} is the horizontal lift of $X$.
Furthermore, by the Chow--Rashevskii theorem \citep[see][Thm.~7.28]{FV10},
\[
G^N(\R^d)=\{\mbX_{0,T}\in T^N(\R^d)\,:\, X \in C^{\var 1}([0,T],\R^d) \}\;,
\]
where again $\mbX_{0,T}$ is defined from $X$ by \eqref{eq:II}.

\begin{definition}[Weakly geometric rough path]
A $p$-rough path $\mbX$ is called \textit{weakly geometric} if $\mbX_{s,t}\in G^N(\R^d)$ for all $(s,t)\in\Delta_T$,
i.e.
\begin{equation}\label{eq:X_shuffle}
\scal{\mbX_{s,t},v\shuffle w} = \scal{\mbX_{s,t},v}\scal{\mbX_{s,t}, w} 
\end{equation}
for all pure tensors $u,v\in T^N(\R^d)$ of combined length at most $N$.
We denote by
$\wgrp_p(V)$
the space of weakly geometric rough paths.
\end{definition}

By the above discussion, $\grp_p(\R^d) \subset \wgrp_p(\R^d)$.
Conversely, by the Chow--Rashevskii theorem and a basic interpolation estimate,
one can approximate every weakly geometric $q$-rough path by smooth paths in the $p$-variation metric for any $p>q\geq 1$,
hence $\wgrp_q(\R^d) \subset \grp_p(\R^d)$.
Viewing $G^N(\R^d)$ as a sub-Riemannian manifold, weakly geometric $p$-rough paths are precisely paths of finite $p$-variations with respect to the Carnot--Carath\'eodory metric,
and within this class, geometric $p$-rough paths admit an analogue of Wiener's characterisation of the closure of smooth paths under $p$-variation norms \citep[see][Sec.~8.6]{FV10}.

If $V$ is infinite-dimensional, the same definition of weakly geometric rough paths can be made provided $G^N(V)$ is defined as the closure of the exponential of $\mfg^N(V)$.
One still has the inclusion $\grp_p(V)\subset \wgrp_p(V)$,
but the converse $\wgrp_q(V)\subset \grp_p(V)$ for $q<p$ is open in general. See \citep{Grong_Nilssen_Schmedding_22_wgrp} where this embedding is shown for $p<3$ and $V$ a Hilbert space.

For $p\geq 2$, the spaces of $p$-rough paths (general, geometric, or weakly geometric) are not vector spaces unless $N=1$ or $V=\R$, essentially due to Chen's identity.

\textbf{Levels $2$ and $3$ and It\^o calculus.}
An important case is when $p\in [2,3)$, which covers the regime of Brownian motion.
In this case, $\mbX$ has two non-trivial components
$\mbX=(1,\mbX^1,\mbX^2)=:(1,X,\bbX)$
and Chen's identity is simply
\begin{align*}
X_{s,t} + X_{t,u} &= X_{s,u}\;,\\
\bbX_{s,t}+\bbX_{t,u}+X_{s,t}\otimes X_{t,u} &= \bbX_{s,u}\;.
\end{align*}
The first condition in particular implies $X_{s,t}=X_{t}-X_s$ for some path $X\colon [0,T]\to \R^d$ that we denote by the same symbol.

Given a Brownian motion
$X\colon[0,T]\to\R^d$,
one can define $\bbX_{s,t}=\int_s^t X_{s,r}\mrd X_{r}$, which is the case $k=2$ of \eqref{eq:II},
and where the integral is taken in either the It\^o or Stratonovich sense.
Then, almost surely, $(1,X,\bbX)$ is a $\gamma$-H\"older rough path for any $\gamma\in (\frac13,\frac12)$ (thus a $\frac{1}\gamma$-rough path).
More generally, any semi-martingale $X$ yields in this way a random $p$-rough path for any $p\in(2,3)$.
Furthermore,
if the Stratonovich integral is used in the definition of $\bbX$, then $\mbX$ is a geometric rough path, while if the It\^o integral is used, $\mbX$ is not geometric.

Two other families of processes that admit stochastic constructions of geometric rough path lifts are Gaussian processes and Markov diffusions.

The former family incorporates Gaussian processes $X\colon [0,T]\to\R^d$ whose covariance function $\rho(s,t) = \E[X_s\otimes X_t]$ satisfies a two-dimensional finite $\rho$-variation condition for $\rho<2$, \citep[see][Chap.~15]{FV10}.
In this case, $X$ admits a canonical geometric $p$-rough path $\mbX$ for $p>2\rho$, whereby canonical means as a limit of natural bounded variation  approximations (e.g. piecewise linear interpolations).
Note that $\rho\geq \frac32$ requires level 3 rough paths.
Fractional Brownian motion with Hurst parameter $H\in(\frac14,\frac12]$
satisfies the required conditions
with $\rho=1/2H$
(the case $H>\frac12$ is covered by Young integration).
For $H\leq \frac14$, there is a fundamental obstruction to defining a canonical geometric $p$-rough path lift of $X$ via bounded variation approximations,
which is the divergence of the variance of the second iterated integrals~\citep[see][]{Coutin_Qian_02}.

The latter family includes
all Markov diffusions with generator $\sum_{i,j=1}^d\partial_i (a^{ij}\partial_j)$ for a measurable, symmetric, uniformly elliptic $a\colon \R^d\to\R^{d\times d}$
(if $a$ lacks sufficient regularity, these diffusions may fail to be semi-martingales).
Through Dirichlet form analysis,
\cite{FV08} showed that such diffusions admit
canonical geometric $\gamma$-H\"older rough path lifts for any $\gamma<\frac12$ (in fact, these lifts can be constructed directly as Markov diffusions associated to Dirichlet forms on $G^N(\R^d)$, which allows for their transition functions to depend on the higher order terms $\mbX^k_{0,t}$, $1\leq k\leq N$).
Markov processes, including diffusions on fractals, were also studied as rough paths by \cite{BHL02}, who show an extension of the classical Wong--Zakai theorem.

\section{Integration and differential equations}

A central result in the theory of rough paths is the following universal limit theorem (ULT).
We denote by $\Psi\colon C^\infty([0,T],V) \hookrightarrow \grp_p(V), X\mapsto \mbX$ the canonical embedding defined by \eqref{eq:II}.
For $\gamma\geq 1$, we denote by $C^\gamma$ the $\gamma$-H\"older  space with bounded $k$-th derivatives for all $k=0,\ldots, \floor \gamma$.

\begin{theorem}[Universal limit theorem]\label{thm:ULT}
Let $V,W$ be Banach spaces. Suppose $1\leq p<\gamma$ and $f\in C^\gamma(W,L(V,W))$.
Then there exists a unique solution to the \emph{rough differential equation (RDE)}
\begin{equation}\label{eq:RDE}
\mrd Y_t = f(Y_t)\mrd \mbX_t\;, \quad Y_0\in W
\end{equation}
for any geometric $p$-rough path $\mbX \in \grp_p(V)$.
The solution map
\[
\hat I\colon \grp_p(V)\times W \to C^{\var p}([0,T],W)\;, \quad \hat I\colon (\mbX,Y_0)\mapsto Y
\]
is locally Lipschitz
and is the unique continuous extension of the classical solution map of \eqref{eq:ODE}  $I\colon (X,Y_0)\mapsto Y$ for smooth $X$
once we embed $C^\infty([0,T],V)$ in $\grp_p(V)$ via $\Psi$.
\end{theorem}

\begin{figure}
\begin{center}
\begin{tikzpicture}[thick,scale=1.5]
\node at ( -1,1.5) [] {$W$};
\node at ( -1,0) [] {$W$};
\node at ( -0.55,0) [] {$\times$};
\node at ( -0.65,1.5) [] {$\times$};
  \node (X) at ( 0,1.5) [] {$\grp_p(V)$};
  \node (O) at ( 0,0) [] {$C^\infty$};
  \node at ( 0,-0.6) [] {$X$};
  \node at ( 0,-0.3) [rotate=90] {$\in$};
  \node at ( 0,1.8) [rotate=-90] {$\in$};
  \node at ( 0,2.1) [] {$\mbX$};

  \node (end1) at ( 3.5,0) [] {$C^{\var p}$};
  \draw [->] (O) -- node [left] {$\Psi$} (X);
  \draw [->] (X) --node [above] {$\hat I$} (end1);
  \draw [->] (O) --node [above] {$I$}  (end1);
\end{tikzpicture}
\caption{Maps from Theorem \ref{thm:ULT}.
Both $\Psi$ and $I$ are discontinuous in $p$-variation norm for $p\geq 2$.}\label{fig:RP_sol_map}
\end{center}
\end{figure}
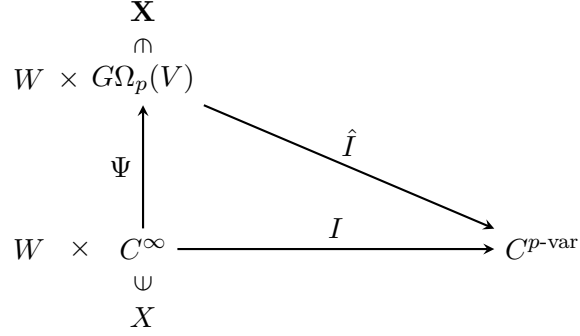

The core of the theorem is the existence of $\hat I$;
its uniqueness is obvious from the density of smooth paths in $ \grp_p(V)$.
See the commutative diagram in Figure \ref{fig:RP_sol_map}.
One should interpret this result as a rough version of the Picard--Lindel\"of theorem.
The solution map $\hat I$ is often called the It\^o--Lyons map.

Several approaches to rough path theory have appeared since its inception by \cite{Lyons98}.
These all share essentially the same definition of a (geometric) rough path,
but they differ in their proof of the ULT
and in their characterisation of the solution $Y$.
We give a description of two of these approaches, the first due to \cite{Gubinelli04}  based on controlled rough paths,
and the second due to \cite{Davie08} based on Euler estimates.

\textbf{Gubinelli's controlled paths.} For a $p$-rough path $X$ and Banach space $W$, a (rough) path controlled by $\mbX$ is a continuous map
\begin{align*}
Y&=(Y^0,\ldots, Y^{N-1})
\colon [0,T]\to \bigoplus_{k=0}^{N-1} L(V^{\otimes k},W)
\end{align*} such that
$\|R^\ell_{s,t}\|\leq \omega(s,t)^{(N-\ell)/p}$ for a control function $\omega$,
where $R^\ell_{s,t}$ is defined, for $0\leq \ell< N$, by
\begin{equation}\label{eq:con_RP_error}
Y_t^\ell = \sum_{k=\ell}^{N-1} Y_{s}^k(\mbX^{k-\ell}_{s,t}) + R^\ell_{s,t}\;.
\end{equation}
Here $Y_{s}^k(\mbX^{k-\ell}_{s,t})\in L(V^{\otimes \ell},W)$ is obtained by $Y_{s}^k(\mbX^{k-\ell}_{s,t})\colon v\mapsto Y^k_s(\mbX^{k-\ell}_{s,t}\otimes v)$.
The components $(Y^1,\ldots, Y^{N-1})$ are called the \textit{Gubinelli derivatives} of $Y$ (or of $Y^0$).
We denote by $D_\mbX^{N/p}$ the space of all paths controlled by $\mbX$, which we equip with the norm
\begin{equation}\label{eq:controlled_norm}
\|Y\|_{\mbX,\var p} = \sum_{\ell=0}^{N-1}\{ \|Y_0^\ell\| +\|R^\ell\|_{\var{p/(N-\ell)};[0,T]}\}\;,
\end{equation}
where $\|R^\ell\|_{\var{q};[0,T]}$ is defined exactly as in \eqref{eq:p_var_linear}.
(The exponent $N/p$ in $D_\mbX^{N/p}$ is made for consistency with \cite{FrizHairer20}
and is motivated by the fact that $Y^0$ behaves as if it had $N$ derivatives with respect to $\mbX$ worth $1/p$-H\"older regularity each.)

\begin{example}\label{ex:CRP}
As a first example, each component of $\mbX$ is a path controlled by $\mbX$.
Indeed,
let $W=V^{\otimes n}$ for $1\leq n\leq N$
and define $Y^k\equiv0$ for $k>n$
and $Y^k_{t} = \mbX^{n-k}_{0,t}$
for $0\leq k\leq n$, which we understand as the linear map $V^{\otimes k}\ni v \mapsto \mbX^{n-k}_{0,t}\otimes v \in V^{\otimes n}$.
Then, by Chen's identity,
\[
\sum_{k=\ell}^{N-1} Y_{s}^k(\mbX^{k-\ell}_{s,t}) 
= \sum_{k=\ell}^{N-1} \mbX^{n-k}_{0,t}\mbX^{k-\ell}_{s,t} = \mbX^{n-\ell}_{0,t}\;,
\]
which is precisely $Y^\ell_t$ unless $n=N$ and $\ell=0$; in this latter case, it is $Y^0_t-\mbX^N_{s,t}$, which still satisfies \eqref{eq:con_RP_error} with $R^0_{s,t}=\mbX^N_{s,t}$.
\end{example}

\begin{example}
As another example, suppose $p\in[1,2)$ so that $N=1$, and $Y^0$ a path of finite $p$-variation.
Then $R^0_{s,t}=Y^0_t-Y^0_s$, which indeed satisfies $\|R^0_{s,t}\|\leq\omega_{Y,p}(s,t)^{1/p}$.
We therefore have equality of Banach spaces $D_\mbX^{1/p} = C^{\var p}$ in the `Young regime' $p\in [1,2)$, .
\end{example}

In contrast to the space of rough paths, $D_\mbX^{N/p}$ is a Banach space; it clearly depends on the underlying rough path $\mbX$.
For fixed $\mbX$, one can define on $D_\mbX^{N/p}$
usual operations of calculus: multiplication, composition by smooth functions, and integration against $\mbX$.
We describe only integration, since this motivates the definition of controlled paths
(concerning multiplication and composition with smooth functions, we only comment that this is where geometricity of $\mbX$ is required in the proof of the ULT).
Suppose $Y^k$ takes values in $L(V^{\otimes( k+1)},W)$ (so the target Banach space of $Y^0$ is $L(V,W)$)
and consider the function
\begin{equation}\label{eq:CRP_integral}
\Xi_{s,t}= \sum_{k=0}^{N-1} Y^k_{s}\mbX^{k+1}_{s,t}\;.
\end{equation}
By Chen's identity
and the definition of $R^\ell$, it follows that $\delta\Xi_{s,u,t}=-\sum_{k=0}^{N-1} R^k_{s,t}\mbX^{k+1}_{t,u}$,
which is of order $\omega(s,u)^{(N+1)/p}$
for a control function $\omega$.
By the sewing lemma, Riemann sums of $\Xi$ over any interval $[0,t]\subset[0,T]$ converge to an object that we denote $\int_0^t Y_s\mrd \mbX_s\in W$.
Furthermore, $\int_0^\bullet Y_s\mrd \mbX_s$ admits a `lift' to a controlled path with Gubinelli derivatives $(Y^0,\ldots, Y^{N-2})$, so $(\int_0^t Y_s\mrd \mbX_s, Y^0,\ldots, Y^{N-2})$ is a path controlled by $\mbX$.

The meaning of the RDE \eqref{eq:RDE} then becomes an integral equation
$Y_t = Y_0 + \int_0^t f(Y_s)\mrd \mbX_s$,
which one can solve for using contraction mapping arguments in the space of paths controlled by $\mbX$.
Under a distance function that generalises the norm \eqref{eq:controlled_norm},
one can further show that the solution map $\mbX\mapsto Y$ is continuous.
If $\mbX$ is the canonical lift of a smooth path $X$, then the $0$-th component $Y^0$ solves the classical ODE $\mrd Y = f(Y)\mrd X$,
from which one sees the commutativity in Figure \ref{fig:RP_sol_map}.
(In fact, one not only enhances the space of `drivers' from paths $X$ to rough paths $\mbX$,
but also the space of solutions from paths $Y^0$ to controlled paths $Y\in D_\mbX^{N/p}$.)
This concludes our sketch proof of the ULT with the approach of Gubinelli -- further details can be found in
\citep{FrizHairer20}.
See in particular \citep[Sec. 8.8]{FrizHairer20}
for comparison with Lyons' original approach based on integration of $1$-forms.

\textbf{Davie's Euler estimates.} Another approach to rough path theory is based on Euler estimates.
We begin with an alternative definition of a solution to the RDE \eqref{eq:RDE}.

\begin{definition}[Davie's RDE solution]
A path $Y\in C^{\var p}([0,T],W)$ is a solution to \eqref{eq:RDE} if there exists a control function $\omega$ and a function $\theta\colon [0,\infty) \to [0,\infty)$ such that $\lim_{\delta\to 0}\theta(\delta) =0$ and
\begin{equation}\label{eq:Euler}
Y_t = Y_s + \sum_{k=1}^{N} f^{\otimes k}(\mbX^k_{s,t})\id (Y_s) + \theta(\omega(s,t))\;,
\end{equation}
where $\id\colon W\to W$ is the identity map, we treat $f$ as an element of $L(V,C^\gamma(W,W))\simeq C^\gamma(W,L(V,W))$, and $f^{\otimes k}$ contracts with $\mbX^k_{s,t}$ to produce an element of $C^\gamma(W,W)^{\otimes k}$
which we treat as a $k$-fold application of the map $(h, g) \mapsto Dg(\cdot)h(\cdot)\colon W \to W$ for $g,h \in C^1(W,W)$.
\end{definition}

For example, for $p\in [3,4)$, \eqref{eq:Euler} becomes
\[
Y_t = Y_s + f(Y_s)\mbX^1_{s,t} + Df(Y_s) f(Y_s) \mbX^2_{s,t} + D[Df \cdot f](Y_s) f(Y_s) \mbX^3_{s,t} + \delta(\omega(s,t))\;.
\]
The estimate \eqref{eq:Euler} is a level-$N$ Euler scheme.
If $X$ is of bounded variation and $\mbX$ is defined by \eqref{eq:II},
one naturally arrives at this scheme by repeated applications of the fundamental theorem of calculus.

Davie's definition is equivalent to stating that $Y$ is a controlled path with Gubinelli derivatives $\{f^{\otimes k}(\cdot)\id (Y_s)\}_{k=1}^{N-1}$,
which satisfies the rough differential equation $Y_t = Y_0 + \int_0^t Y_s \mrd \mbX_s$.
It follows that Davie's definition of solution is equivalent to the one previously seen with Gubinelli's approach (but does not require the notion of a rough integral in its formulation).

Working in finite dimensions $V=\R^d$, we give a brief sketch of Davie's approach to showing existence of solutions to RDEs.
The next two estimates are central in this regard.
Suppose $x,\tilde x$ are bounded variation paths with canonical lifts $\mbx,\tilde\mbx\in \grp_p(\R^d)$
such that $\|\tilde x\|_{\var 1;[s,t]}\leq K \|\mbx\|_{\var p;[s,t]}$ for some $K\ge 1$ and $\tilde\mbx_{s,t} = \mbx_{s,t}$.
Let $y,\tilde y$ be solutions to the ODE \eqref{eq:ODE} driven by $x,\tilde x$ respectively over the interval  
$[s,t]$ with the same initial value $y_s=\tilde y_s$.
Then
\begin{equation}\label{eq:a-priori}
\|y\|_{\var p;[s,t]} \lesssim \phi_p(\|f\|_{C^{\gamma-1}}\|\mbx\|_{\var p;[s,t]})
\end{equation}
where $\phi_p(z)=z+z^p$,
and
\begin{equation}\label{eq:y_diff}
|\tilde y_t - y_t| \lesssim (K\|f\|_{C^{\gamma-1}}\|\mbx\|_{\var p;[s,t]})^\gamma\;.
\end{equation}
These two estimates are known as Davie's lemma.
Similar to the previously discussed approach based on the sewing lemma, an important step in the proof of \eqref{eq:a-priori}-\eqref{eq:y_diff} is to bound the quantity $\delta\Xi_{u,v,w}$ where $\Xi_{u,v} = y_{u,v}-y^{(u,v)}_{u,v}$
and $y^{(u,v)}$ is the solution to \eqref{eq:ODE} over $[u,v]$ with initial value $y^{(u,v)}_u= y_u$ and driven by $x^{(u,v)}$ that satisfies the same properties as $\tilde x$ above for the interval $[u,v]$;
such $x^{(u,v)}$ always exists by taking the sub-Riemannian geodesic in $G^N(\R^d)$ connecting $1$ and $\mbx_{u,v}$
\citep[see][Chap.~8]{FV10}.

With the estimates \eqref{eq:a-priori}-\eqref{eq:y_diff} in hand, one can construct solutions in the sense of \eqref{eq:Euler}
for any $\mbX\in \wgrp_p(\R^d)$ by approximating $\mbX$ with canonical lifts $\mbx^{(n)}$ of bounded variation paths $x^{(n)}$ such that $\mbx^{(n)}\to\mbX$ in the uniform norm and $\sup_n\|\mbx^{(n)}\|_{\var p;[0,T]}<\infty$.
The bound \eqref{eq:a-priori} gives equicontinuity of the solutions $y^{(n)}$ to \eqref{eq:ODE} driven by $x^{(n)}$, thus convergence $y^{(n)}\to Y$ along a subsequence to a candidate solution $Y$.
By taking $\tilde x$ in \eqref{eq:y_diff} as, for example, the geodesic connecting $1\in G^N(\R^d)$ and $\mbx^{(n)}_{s,t}$,
one can verify that the candidate solution $Y$ indeed satisfies \eqref{eq:Euler}.

An advantage of Davie's approach is that it provides \textit{existence} of solutions in the regime $f\in C^{\gamma-1}$ for $\gamma>p$ (instead of $f\in C^\gamma$)
wherein one does not in general have uniqueness
(akin to the classical Peano existence theorem vs. Picard--Lindel\"of theorem).
With the extra regularity $f\in C^\gamma$ in Theorem \ref{thm:ULT},
one can refine the approach to show uniqueness together with the claimed local Lipschitz continuity.

Given a solution $Y$ to an RDE (in any sense discussed above and for any Banach spaces $V,W$),
one always has an Euler-type estimate for $Y_t$ as $t\downarrow0$ that is of order $(\|f\|_{C^{\gamma-1}}\|\mbX\|_{\var p;[0,t]})^{\gamma}$ and which uses only $\mbX_{0,t}$ (or its level-$\floor\gamma$ extension $\hat \mbX_{0,t}$ if $\gamma\geq \floor p +1$ - see below).

\textbf{Optimality.}
We comment on several ways that the ULT may be optimised and extended.

\begin{itemize}
\item If $f$ is linear (thus smooth but not bounded),
global solutions to \eqref{eq:RDE} still exist
with precise estimates on the $p$-variation of $Y$.

\item The H\"older-regularity assumption $f\in C^\gamma$ may also be weakened to Lipschitz-regularity $f\in \Lip^\gamma$ in the sense of E. Stein;
for $\gamma$ non-integer, $\Lip^\gamma$ agrees with $C^\gamma$, but for $\gamma$ an integer, it only requires that the $k$-th derivative for $k=0,\ldots,\gamma-1$ is Lipschitz and bounded.

\item If $\gamma=p$, one can still has existence and uniqueness of solutions $Y$ to \eqref{eq:RDE}
together with continuity of $\mbX\mapsto Y$ (at least for $V=\R^d$), although one loses local Lipschitz regularity.
In fact, this result holds for $f\in \Lip^p$ and where we lower the regularity assumptions on $\mbX$ from finiteness of $p$-variation in \eqref{eq:p-var} to finiteness of
\begin{equation*}
\omega_{\mbX,\psi_p}(s,t) = \sup_{P}  \sum_{[u,v]\in P} \psi_{p}
\Big(
\max_{1\leq k\leq N} \|\mbX^k_{u,v}\|_{V^{\otimes k}}^{1/k}
\Big)\;,
\end{equation*}
where $\psi_p(z)=z^p/(\log |\log z|)$ for $z\ll 1$.
This is of some significance because, if $\mbX$ is the canonical lift of Brownian motion, then $\omega_{\mbX,\psi_2}(s,t)<\infty$ but $\mbX$ has infinite $2$-variation almost surely.
\end{itemize}

\textbf{Applications.}
One of the first applications of rough path theory in stochastic analysis is to give meaning to differential equations driven by irregular signals that are not semi-martingales and thus fall outside the scope of It\^o calculus and Young integration.
This includes fractional Brownian motion with Hurst parameter $H\in (\frac14,\frac12)$ \citep[see][]{Coutin_Qian_02} and symmetric Markov diffusions \citep[see][]{BHL02,FV08}.
Due to the ULT, these solution theories automatically come with pathwise approximations results.
In the Gaussian case, the ULT furthermore intertwines elegantly with Malliavin calculus: \cite{CF10} showed existence of densities to RDEs driven by Gaussian rough paths, thus giving rise to non-Markovian H\"ormander theory,
and this was later improved to smoothness of densities by \cite{CHLT15}.

It turns out that, under natural assumptions, the rough integral and solutions to rough differential equations coincide almost surely with their counterparts in It\^o/Stratonovich calculus when the rough path $\mbX$ is taken as the canonical enhancement of a semi-martingale.
This link, together with the ULT, allows for a number of simplified
proofs of results in stochastic analysis of semi-martingales, e.g. Wong-Zakai convergence theorem, Stroock--Varadhan support theorem, Freidlin--Wentzell large deviations principles \citep[see][Chapters~17,~19]{FV10}.
The simplification that rough paths bring is to factor the solution map as $I = \hat I\circ \Psi$ where $\hat I$ is continuous by the ULT and $\Psi$ is not continuous but is simpler to analyse than the full solution map $I$ (indeed, $\Psi(X)$ is the solution to a single linear SDE driven by $X$).
Rough path theory further allows for new pathwise approaches to McKean--Vlasov equations \citep[see][]{Cass_Lyons_15,CDFM20}
and to stochastic partial differential equations (SPDEs)
driven by irregular noise,
see \citep[Sec. 12]{FrizHairer20} and also below for generalisations to \textit{singular} SPDEs.

Rough path theory was applied by \cite{KM16,KM17} in the study of deterministic fast-slow systems of the form $\dot x = a(x,y) + \eps^{-1}b(x,y), \dot y = \eps^{-2}g(y)$ to show that the slow component $x$ converges in law to a solution of an explicit stochastic differential equation provided that $g$ generates a suitable ergodic flow with invariant probability measure $\mu$ and $y_0$ is sampled according to $\mu$.
The work of \cite{KM17} in particular relies on rough path theory in infinite-dimensional spaces.
See also \citep{BCKL19} where averaging of electron speeds is shown for a variant of the Drude--Lorentz model by means of rough paths.

See \citep{FrizHairer20} and references therein for further stochastic analytic applications of rough path theory.

\section{Extensions}

Given a $p$-rough path $\mbX$ and $M\geq N=\floor p$, there is a unique extension of $\mbX$ to a map $\hat \mbX \colon\Delta_T \to T^M(V)$,
which is called the (level-$M$) \textit{multiplicative extension} of $\mbX$ by \cite{Lyons98},
that satisfies \ref{pt:Chen} (Chen's identity)
and \ref{pt:pvar} (with $N$ replaced by $M$ and a possibly larger $\omega$).
By `extension', we mean $\hat \mbX^k=\mbX^k$ for all $k\leq N$.
The extension $\hat \mbX$ furthermore satisfies the following estimate: for all $\beta>0$ sufficiently large, if $\omega$ is a control function for which
\begin{equation}\label{eq:p_var_sharp}
\max_{1\leq k\leq N}(\beta(k/p)!\|\mbX^k_{s,t}\|_{V^{\otimes k}})^{p/k} \leq \omega(s,t)
\end{equation}
for all $(s,t)\in\Delta_T$,
then \eqref{eq:p_var_sharp} also holds with $\mbX$ replaced by $\hat \mbX$ and $N$ replaced by $M$.
Above, $x!=\Gamma(x+1)$ for the classical gamma function $\Gamma$.
\cite{Hara_Hino_10} showed that
a lower bound for $\beta$ is $\beta\geq p(1+2^{\gamma}(\zeta(\gamma)-1))$, where $\gamma=\frac{\floor p + 1}{p} > 1$ and $\zeta(s) = \sum_{n=1}^\infty n^{-s}$ is the zeta function
(\cite{Lyons98} proved the slightly weaker bound $\beta\geq p^2(1+2^{\gamma}(\zeta(\gamma)-1))$).
The extension map $\mbX\mapsto \hat \mbX$ is further continuous in the $p$-variation metric \eqref{eq:d_p-var}.

Note that \eqref{eq:p-var-omega} and \eqref{eq:p_var_sharp} are essentially equivalent in that one implies the other after possibly increasing $\omega$ by a constant (depending on $p$).
In particular, every $p$-rough path is also canonically a $q$-rough path for $q\geq p$.
By taking smooth approximations, it is immediate that the extension of a geometric $p$-rough path is again geometric.
Likewise the extension of a weakly geometric rough path is again weakly geometric - in finite dimensions, this is clear be seen from the inclusion $\wgrp_q(\R^d)\subset \grp_p(\R^d)$ for $p>q$, while in infinite dimensions, this was shown by \cite{Cass_Driver_Lim_16_wgrp}.

Perhaps most simply, this extension is obtained as the solution of the linear rough differential equation
\begin{equation}\label{eq:RDE_extension}
\mrd \hat \mbX = \hat \mbX \otimes \mrd \mbX\;, \quad \hat \mbX_0 = 1 \in T^{M}(V)\;,
\end{equation}
which exists and is unique.
The extension can also be constructed inductively, where the step from  level $N$ to level $N+1$ is given by integration of controlled rough paths as follows: taking the controlled path from Example \ref{ex:CRP} with $n=N$, we can define $\hat \mbX^{N+1}_{0,t}$ by
`sewing' the two-parameter map
$\Xi_{r,s} = \sum_{k=0}^{N-1} \mbX^{N-k}_{r}\mbX^{k+1}_{r,s} \in V^{\otimes (N+1)}$,
to which the sewing lemma applies since it is precisely of the form \eqref{eq:CRP_integral}.

In addition to the extension problem of a $p$-rough path, it is natural to ask whether a $V$-valued path
admits a lift to a $p$-rough path.
An answer was given by \cite{Lyons_Victoir_07} who showed that, for $p\geq 1$ such that $p \notin \Z\setminus\{1\}$,
if $X\colon [0,T]\to V$ is of finite $p$-variation, then there exists a weakly geometric $p$-rough path lift of $X$ (their result further covers partial extensions, i.e. when the components $\mbX^j$, $j=1,\ldots, k$ for $k<\floor p$ are given with the natural regularity and one is asked to extend $\mbX$ to a weakly geometric $p$-rough path).
If $p>1$, the assumption that $p$ is non-integer cannot be dropped  even when $V=\R^d$ unless $d=1$ \cite[see][Sec. 9.2]{FV10}.

\section{Signature}

The \textit{signature} of a $p$-rough path $\mbX$ is the sequence of tensors
\[
S(\mbX) := (1, S^1(\mbX),S^2(\mbX),\ldots) := (1,\mbX^1_{0,T},\mbX^2_{0,T}\ldots) \in T((V)) := \prod_{k=0}^\infty V^{\otimes k}\;,
\]
where $\mbX^k_{0,T}\in V^{\otimes k}$ for every $k>N$
is obtained by extending $\mbX$ to some arbitrary order $M\geq k$.

Denoting $X_t=\mbX^1_{0,t}$, if $p\in [1,2)$,
the signature $S(X)\equiv S(\mbX)$ is completely determined by $X$ and is given by the sequence of iterated integrals \eqref{eq:II}.

For $X$ smooth,
the study of the map $X\mapsto S(X)$ has a history dating back to Chen, Ree, Magnus in the 50's.
We discuss several algebraic, geometric, and probabilistic properties of the signature below.

\textbf{Algebraic properties.} Supposing for now that $V=\R^d$ is finite dimensional and $X$ is smooth, a fundamental property of the signature is that $\log S(X)$, where $\log$ is taken in $T((\R^d))$ via a formal power series,
takes values in the space of Lie series generated by $\R^d$;
this result is a generalisation of the Baker--Campbell--Hausdorff formula due to Chen.
This turns out equivalent to the fact discussed above that $S(X)$ satisfies the shuffle identities
\begin{equation}\label{eq:shuffles}
\scal{S(X),u\shuffle v}=\scal{S(X),u}\scal{S(X),v}\;,
\end{equation}
i.e. $S(X)$ belongs to $G(\R^d)$ defined by \eqref{eq:G_def}.
By density, signatures of geometric rough paths also take values in $G(\R^d)$.

\textbf{Geometric and analytic properties.}
At the geometric level,
Chen showed that $S(X)$ determines the underlying path $X$, up to translation and reparametrisation, within a class of piecewise differentiable paths.
Using hyperbolic development, \cite{HL10} extended this result to bounded variation paths $X,Y\colon[0,T]\to\R^d$, showing that $S(X)=S(Y)$ if and only if $X$ and $Y$ are tree-like equivalent.
To define tree-like equivalence, we say that a continuous path $Z\colon[0,T]\to\R^d$ is tree-like if there exists a real tree $\mfT$ \citep[see][]{Bestvina_02_trees}, a continuous map $\phi\colon [0,T]\to\mfT$, and a map $\psi\colon \mfT\to\R^d$ such that $Z=\psi\circ\phi$.
We then say that $X,Y\colon[0,T]\to \R^d$ are tree-like equivalent if $X*\overleftarrow{Y}$ is tree-like, where
$X*\overleftarrow{Y}\colon [0,2T]\to\R^d$ is the concatenation of $X$ with the reversal of $Y$,
i.e. $(X*\overleftarrow{Y})_t = X_t$ if $t\in [0,T]$ and $(X*\overleftarrow{Y})_t = X_T-Y_T+Y_{2T-t}$ if $t\in (T,2T]$.

Using different techniques,
\cite{BGLY16}
extended the equivalence $S(X)=S(Y) \Leftrightarrow$ $X$ and $Y$ are tree-like equivalent to all geometric $p$-rough paths (this final result applies to any Banach space $V$, not necessarily finite-dimensional).

The local and asymptotic behaviour of the series $S^k(X)$ furthermore carries analytic information of the underlying path.
If $X\colon [0,1]\to\R^d$ is $C^1$, parametrised at constant speed $|\dot X| \equiv l$,
then \cite{HL10} showed that there exists $N>0$, depending on $l$ and  the modulus of continuity $\delta_{\dot X}$ of $\dot X$,
such that $S^k(X)\neq 0$ for some $1\leq k \leq N$.
They further show that if $\lim_{\eps\to 0}\eps^{-3/4}\delta_{\dot X}(\eps) = 0$,
then
\[
\lim_{k\to\infty}( k! \|S^k(X)\|_{(\R^d)^{\otimes k}} )^{1/k}= \|X\|_{\var 1}\;.
\]
Extending this result to all bounded variation paths currently remains open.

\cite{LyonsXu18} have devised an algorithm based on symmetrisation to recover a $C^1$ path from its signatures.
More precisely, given a $C^1$ path $X$  parametrised at constant speed for the $\ell^1$ norm on $\R^d$, they construct a piecewise linear path $X^{(N)}$ using the truncated signature $(1,S^1(X),\ldots,S^N(X))$ such that $\|X^{(N)}-X\|_{\Hol1}\to 0$ as $N\to\infty$.

\textbf{Probabilistic properties.}
At the probabilistic level, the expectation of the signature for a random path is a natural generalisation of moments.
Indeed, for $x\in\R^d$ and the straight-line path $X\colon[0,1]\to \R^d, t\mapsto tx$ connecting $0$ and $x$,
the signature is
\[
S(X)=\exp(x) = (1,x,x^{\otimes 2}/2!,x^{\otimes 3}/3!,\ldots )\;.
\]
In this way, and in light of the shuffle identities \eqref{eq:shuffles},
the signature $S(X)$ for a path of bounded variation $X\colon[0,T]\to\R^d$ provides a generalisation of polynomials from vectors to paths.
If $x$ is an $\R^d$-valued random variable with moments of all orders,
and $X_t= tx$ for $t\in[0,1]$,
then the (component-wise) \textit{expected signature}
$
\E S(X)
$
encodes the moments of $x$.
Therefore, under suitable integrability conditions (e.g. the essentially sharp Carleman's condition on the rate of growth of $\E[x^{\otimes k}]$),
$\E S(X)$ determines the law of $x$.
There is a generalisation of this fact to any random geometric $p$-rough path $\mbX$ for $p\in [1,\infty)$:
if the expected signature
\[
\E S(\mbX) = (1,\E [\mbX^1_{0,T}],\E [\mbX^2_{0,T}],\ldots) \in T((\R^d))\;,
\]
exists and has an infinite radius of convergence, by which we mean
\begin{equation}\label{eq:tails}
\forall \lambda>0\,:\quad
\sum_{k=0}^\infty \lambda^k \|\E [\mbX^k_{0,T}]\|_{(\R^d)^{\otimes k}} <\infty\;,
\end{equation}
then $\E S(\mbX)$ uniquely determines the law of the random signature $S(\mbX)$ (and thus the law of $\mbX$ up to tree-like equivalence).

This result provides an answer to the uniqueness part of the `moment problem' for random signatures
(\eqref{eq:tails} is a sufficient condition for $\E S(\mbX)$ to determine the law of $S(\mbX)$, but it is far from necessary, which is clear for $d=1$ by comparing it with Carleman's condition).
The result was shown by \cite{CL16} using a non-commutative Fourier transform that is based on 
(Cartan) developments of paths into unitary Lie groups and relies on (non-existence of) Lie polynomial identities in unitary Lie algebras.
Like in the classical case, this Fourier transform always determines the law of a random signature
(without any integrability assumptions).
An earlier solution to the `moment problem' was given by Fawcett, essentially in the case that the law of $\mbX$ has compact support, in which case determinacy of the law from $\E S(\mbX)$ follows from the Stone--Weierstrass theorem and the shuffle identities \eqref{eq:shuffles}.

The expected signature (and non-commutative Fourier transform) can be computed explicitly for
L\'evy processes suitably interpreted as geometric rough paths \citep[see][]{FS17,Chevyrev18},
which gives a direct way to verify condition \eqref{eq:tails};
it in particular holds if $\mbX$ is the canonical lift of Brownian motion.
Condition \eqref{eq:tails} can also be verified for classes of Gaussian and Markovian rough paths, but is known to fail for the signature of Brownian motion stopped upon exiting a smooth domain \citep[see][]{Li_Ni_22}.

Due to its role as the generalisation of polynomials to pathspace, the signature has received interest as a feature of time-ordered data in machine learning.
For example, a maximum mean discrepancy (a form of metric) on probability measures on geometric $p$-rough paths was proposed by \cite{CO22} based on inner products of normalised signatures.
In turn, \cite{KO19} showed that inner products of (normalised) signatures can be computed efficiently using dynamic programming and have applied this to classification problems.
See \citep{CK16, FLMS_23} for further applications of signatures and rough paths in machine learning.

\section{Generalisations}

There have been a number of generalisations of rough path theory in different directions.
We focus only on a few of these -- further discussions can be found in \citep{FV10,FrizHairer20}.

\textbf{Branched rough paths.}
An important generalisation of geometric rough paths that allows for a version of the universal limit theorem
while relaxing the geometricity condition \eqref{eq:X_shuffle}
is branched rough paths of \cite{Gubinelli_10}.
The starting point for this construction is to switch the target space from a truncated tensor algebra to a truncated Grossman--Larson algebra of labelled trees (dual of the Connes--Kreimer algebra).
We give a description of this algebra for $V=\R^d$,
although an infinite dimensional setting is also possible.

A \textit{labelled tree} is a tree $\tau$ (connected graph without cycles) with vertex set $V$ together with a map $\mfn \colon V\to [d]:=\{1,\ldots, d\}$
and a distinguished vertex $\rho_\tau\in V$ called the root.
When we draw trees, we write the label $\mfn(v)$ next to each vertex $v$.
For example, the trees with a single vertex are
\[
\bullet_1,\ldots, \bullet_d\;,
\]
and the trees with two and three vertices are, respectively,
\[
\begin{tikzpicture}[scale=0.2,baseline=0.1cm]
        \node at (0,0)  [dot,label= {[label distance=-0.2em]below: \scriptsize  $ i $} ] (root) {};
         \node at (0,2)  [dot,label={[label distance=-0.2em]above: \scriptsize  $ j$}] (right) {};
            \draw[kernel1] (right) to
     node [sloped,below] {\small }     (root);
     \end{tikzpicture}
\;,
\quad
i,j \in [d]\;,\qquad
\begin{tikzpicture}[scale=0.2,baseline=0.1cm]
        \node at (0,0)  [dot,label= {[label distance=-0.2em]below: \scriptsize  $ i $} ] (root) {};
         \node at (1,2)  [dot,label={[label distance=-0.2em]above: \scriptsize  $ k $}] (right) {};
         \node at (-1,2)  [dot,label={[label distance=-0.2em]above: \scriptsize  $ j $} ] (left) {};
            \draw[kernel1] (right) to
     node [sloped,below] {\small }     (root); \draw[kernel1] (left) to
     node [sloped,below] {\small }     (root);
     \end{tikzpicture}\;,
\quad i,j,k \in [d]\;.
\]
We identify trees up to root and label-preserving isomorphisms; in particular, the order of branches drawn above does not matter, i.e.
$\begin{tikzpicture}[scale=0.2,baseline=0.1cm]
        \node at (0,0)  [dot,label= {[label distance=-0.2em]below: \scriptsize  $ i $} ] (root) {};
         \node at (1,2)  [dot,label={[label distance=-0.2em]above: \scriptsize  $ k $}] (right) {};
         \node at (-1,2)  [dot,label={[label distance=-0.2em]above: \scriptsize  $ j $} ] (left) {};
            \draw[kernel1] (right) to
     node [sloped,below] {\small }     (root); \draw[kernel1] (left) to
     node [sloped,below] {\small }     (root);
     \end{tikzpicture}
=
\begin{tikzpicture}[scale=0.2,baseline=0.1cm]
        \node at (0,0)  [dot,label= {[label distance=-0.2em]below: \scriptsize  $ i $} ] (root) {};
         \node at (1,2)  [dot,label={[label distance=-0.2em]above: \scriptsize  $ j $}] (right) {};
         \node at (-1,2)  [dot,label={[label distance=-0.2em]above: \scriptsize  $ k $} ] (left) {};
            \draw[kernel1] (right) to
     node [sloped,below] {\small }     (root); \draw[kernel1] (left) to
     node [sloped,below] {\small }     (root);
     \end{tikzpicture}
$.

A \textit{forest} is an unordered collection of labelled trees $\tau_1 \cdots \tau_n$, including the empty forest denoted by $1$ that corresponds to $n=0$.
Every labelled tree $\tau$ admits a unique representation $\tau=[\tau_1 \cdots \tau_n]_i$, where $i=\mfn(\rho_\tau)$, $n$ is the degree of $\rho_\tau$, and
$\tau_1 \cdots \tau_n$ is a forest in which each
$\tau_i$ is a labelled tree whose root is adjacent to $\rho_\tau$
(this inductive structure can be used as an alternative definition of a tree that does not make reference to graphs).
For example, 
$\bullet_i=[1]_i$
and
$\begin{tikzpicture}[scale=0.2,baseline=0.1cm]
        \node at (0,0)  [dot,label= {[label distance=-0.2em]below: \scriptsize  $ i $} ] (root) {};
         \node at (1,2)  [dot,label={[label distance=-0.2em]above: \scriptsize  $ k $}] (right) {};
         \node at (-1,2)  [dot,label={[label distance=-0.2em]above: \scriptsize  $ j $} ] (left) {};
            \draw[kernel1] (right) to
     node [sloped,below] {\small }     (root); \draw[kernel1] (left) to
     node [sloped,below] {\small }     (root);
     \end{tikzpicture}
=
[\bullet_j \bullet_k]_i
$.

The Connes--Kreimer (CK) Hopf algebra $H_\CK$ is the formal vector space spanned by the set of forests equipped
with the (commutative) forest product, defined for $\tau=\tau_1\dots\tau_n$ and $\sigma=\sigma_1\cdots\sigma_m$ by
\[
\tau\sigma  =
\tau_1\cdots\tau_n \sigma_1 \cdots\sigma_m\;,
\]
for which the empty forest $1$ is the identity,
and with the Connes--Kreimer coproduct $\Delta_\CK$,
which is defined
inductively by $\Delta_\CK 1 = 1\otimes 1$ and for every $\tau\in H_\CK$ and $i\in [d]$ by
\begin{equation}\label{eq:CK}
\Delta_\CK [\tau]_i = (\id\otimes [\cdot]_i)(\Delta_\CK \tau) + [\tau]_i\otimes 1\;,
\end{equation}
where we extend the map $[\cdot]_i \colon \tau \mapsto [\tau]_i$ to $H_\CK$ by linearity,
together with $\Delta_\CK(\tau_1\cdots\tau_n) = (\Delta_\CK\tau_1)\cdots(\Delta_\CK\tau_n)$.
Concretely, 
for a tree $\tau$,
\[
\Delta_\CK \tau =\tau\otimes 1 + 1\otimes\tau+ \sum_{c} \tau_c\otimes  \tau/c\;,
\]
where the sum is over all admissible cuts of $\tau$.
Here, an admissible cut $c$ is a non-empty set of edges of $\tau$ such that, for any vertex $x$ of $\tau$,
the unique path connecting $x$ to the root of $\tau$ traverses at most one edge in $c$.
Then $\tau_c$ is the forest obtained by collecting all the subtrees of $\tau$ that are above an edge in $c$, and $\tau/c$ is the tree that remains after all subtrees in $\tau_c$ are removed, including the edges in $c$ (we prune the branches connecting trees in $\tau_c$ to the trunk $\tau/c$). For example,
\begin{align*}
\Delta_\CK
\begin{tikzpicture}[scale=0.2,baseline=0.1cm]
        \node at (0,0)  [dot,label= {[label distance=-0.2em]below: \scriptsize  $ i $} ] (root) {};
         \node at (1,2)  [dot,label={[label distance=-0.2em]above: \scriptsize  $ k $}] (right) {};
         \node at (-1,2)  [dot,label={[label distance=-0.2em]left: \scriptsize  $ j $} ] (left) {};
\node at (-1,4)  [dot,label={[label distance=-0.2em]left: \scriptsize  $ \ell $} ] (up) {};
            \draw[kernel1] (right) to
     node [sloped,below] {\small }     (root);
\draw[kernel1] (left) to
     node [sloped,below] {\small }     (root);
\draw[kernel1] (up) to
     node [sloped,below] {\small }     (left);
     \end{tikzpicture}
=&
\begin{tikzpicture}[scale=0.2,baseline=0.1cm]
        \node at (0,0)  [dot,label= {[label distance=-0.2em]below: \scriptsize  $ i $} ] (root) {};
         \node at (1,2)  [dot,label={[label distance=-0.2em]above: \scriptsize  $ k $}] (right) {};
         \node at (-1,2)  [dot,label={[label distance=-0.2em]left: \scriptsize  $ j $} ] (left) {};
\node at (-1,4)  [dot,label={[label distance=-0.2em]left: \scriptsize  $ \ell $} ] (up) {};
            \draw[kernel1] (right) to
     node [sloped,below] {\small }     (root);
\draw[kernel1] (left) to
     node [sloped,below] {\small }     (root);
\draw[kernel1] (up) to
     node [sloped,below] {\small }     (left);
     \end{tikzpicture}
\otimes 1
+
1\otimes
\begin{tikzpicture}[scale=0.2,baseline=0.1cm]
        \node at (0,0)  [dot,label= {[label distance=-0.2em]below: \scriptsize  $ i $} ] (root) {};
         \node at (1,2)  [dot,label={[label distance=-0.2em]above: \scriptsize  $ k $}] (right) {};
         \node at (-1,2)  [dot,label={[label distance=-0.2em]left: \scriptsize  $ j $} ] (left) {};
\node at (-1,4)  [dot,label={[label distance=-0.2em]left: \scriptsize  $ \ell $} ] (up) {};
            \draw[kernel1] (right) to
     node [sloped,below] {\small }     (root);
\draw[kernel1] (left) to
     node [sloped,below] {\small }     (root);
\draw[kernel1] (up) to
     node [sloped,below] {\small }     (left);
     \end{tikzpicture}
+
\bullet_\ell
\otimes
\begin{tikzpicture}[scale=0.2,baseline=0.1cm]
        \node at (0,0)  [dot,label= {[label distance=-0.2em]below: \scriptsize  $ i $} ] (root) {};
         \node at (1,2)  [dot,label={[label distance=-0.2em]above: \scriptsize  $ k $}] (right) {};
         \node at (-1,2)  [dot,label={[label distance=-0.2em]above: \scriptsize  $ j $} ] (left) {};
            \draw[kernel1] (right) to
     node [sloped,below] {\small }     (root);
\draw[kernel1] (left) to
     node [sloped,below] {\small }     (root);
\end{tikzpicture}
+
\begin{tikzpicture}[scale=0.2,baseline=0.1cm]
        \node at (0,0)  [dot,label= {[label distance=-0.2em]left: \scriptsize  $ j $} ] (root) {};
         \node at (0,2)  [dot,label={[label distance=-0.2em]left: \scriptsize  $ \ell$}] (right) {};
            \draw[kernel1] (right) to
     node [sloped,below] {\small }     (root);
     \end{tikzpicture}
\otimes
\begin{tikzpicture}[scale=0.2,baseline=0.1cm]
        \node at (0,0)  [dot,label= {[label distance=-0.2em]left: \scriptsize  $ i $} ] (root) {};
         \node at (0,2)  [dot,label={[label distance=-0.2em]left: \scriptsize  $ k$}] (right) {};
            \draw[kernel1] (right) to
     node [sloped,below] {\small }     (root);
     \end{tikzpicture}
\\
&+
\bullet_k
\otimes
\begin{tikzpicture}[scale=0.2,baseline=0.1cm]
        \node at (0,0)  [dot,label= {[label distance=-0.2em]left: \scriptsize  $ i $} ] (root) {};
         \node at (0,2)  [dot,label={[label distance=-0.2em]left: \scriptsize  $ j$}] (right) {};
\node at (0,4)  [dot,label={[label distance=-0.2em]left: \scriptsize  $ \ell$}] (up) {};
            \draw[kernel1] (right) to
     node [sloped,below] {\small }     (root);
 \draw[kernel1] (up) to
     node [sloped,below] {\small }     (right);
     \end{tikzpicture}
+
\bullet_\ell\bullet_k
\otimes
\begin{tikzpicture}[scale=0.2,baseline=0.1cm]
        \node at (0,0)  [dot,label= {[label distance=-0.2em]below: \scriptsize  $ i $} ] (root) {};
         \node at (0,2)  [dot,label={[label distance=-0.2em]above: \scriptsize  $ j$}] (right) {};
            \draw[kernel1] (right) to
     node [sloped,below] {\small }     (root);
     \end{tikzpicture}
+
\begin{tikzpicture}[scale=0.2,baseline=0.1cm]
        \node at (0,0)  [dot,label= {[label distance=-0.2em]below: \scriptsize  $ j $} ] (root) {};
         \node at (0,2)  [dot,label={[label distance=-0.2em]above: \scriptsize  $ \ell$}] (right) {};
            \draw[kernel1] (right) to
     node [sloped,below] {\small }     (root);
     \end{tikzpicture}
\bullet_k\otimes
\bullet_i
\;.
\end{align*}
Finally, the Grossman--Larson (GL) algebra $H_\GL$ is the dual of the CK Hopf algebra, i.e. $H_\GL=H_\CK^*$ as a vector space
equipped with the GL product $\star$ defined by $\scal{\tau^*\star\sigma^*,\tau} = \scal{\tau^*\otimes\sigma^*,\Delta_\CK \tau}$ for all $\tau^*,\sigma^*\in H_\GL$ and $\tau \in H_\CK$.
As a vector space, $H_\GL$ is isomorphic to the vector space of formal series of forests with real coefficients.

For a path of bounded $p$-variation $X\in C^{\var p}([0,T],\R^d)$ with $p\in[1,2)$ and $(s,t)\in\Delta_T$,
there is a natural way to assign a element $\mbX_{s,t}\in H_\GL$
by
\begin{equation}\label{eq:mult_CK}
\scal{\mbX_{s,t},\tau\sigma}=\scal{\mbX_{s,t},\tau}\scal{\mbX_{s,t},\sigma}
\end{equation}
for all forests $\tau,\sigma$ (so that $\mbX_{s,t}$ is a character on $H_\CK$)
and
inductively for trees by
\begin{equation}\label{eq:II_branched}
\scal{\mbX_{s,t},[\tau]_i} = \int_s^t  \scal{\mbX_{s,r},\tau} 
\mrd X^i_r\;,
\end{equation}
where the integral is in the Young sense.
Identity \eqref{eq:mult_CK} is the analogue of the shuffle relations \eqref{eq:shuffles}
(although \eqref{eq:mult_CK} is a definition rather than a non-trivial identity).
Furthermore, it readily follows by induction from the definition of $\Delta_\CK$ \eqref{eq:CK} that
\[
\mbX_{s,t} \star \mbX_{t,u} = \mbX_{s,u}\;,
\]
which is the analogue of Chen's identity \ref{pt:Chen} in Definition~\ref{def:RP}.

Let $|\tau|$ denote the number of vertices in a tree $\tau$, and likewise $|\tau_1\cdots\tau_n|=|\tau_1|+\cdots+|\tau_n|$,
and $F_N$ denote the set of forests $\tau$ with $|\tau|\leq N$.
For $N\geq 1$ an integer, we let $J_N$ denote the set of all $x\in H_\GL$ such that $\scal{x,\tau}=0$ for all $\tau\in F_N$.
Then $J_N$ is an ideal in $H_\GL$ and we let $H^N_{\GL} = H_\GL/J_N$ denote the truncated Grossman--Larson algebra,
which is canonically isomorphic (as a vector space) to the
dual of $\Span_\R\{\tau\,:\,\tau\in F_N\}$.
We define
\[
B^N = \{g\in H^N_{\GL}\,:\,\scal{g,1}=1\,,\, \scal{g,\tau\sigma}=\scal{g,\tau}\scal{g,\sigma}
 \,\forall \tau,\sigma\in F_N \text{ such that } \tau\sigma\in F_N\}\;.
\]
Then $B^N$ is a group with multiplication $\star$ and is the truncated version of the Butcher group (the group of characters on $H_\CK$) that arises in numerical analysis.

With these preliminaries in place, for $p\geq 1$,
a branched $p$-rough path is a continuous map
$\mbX \colon \Delta_T \to B^N$, where $N=\floor{p}$ as usual,
which satisfies
\begin{enumerate}[label=(\alph*B)]
\item\label{pt:Chen_branched} Chen's identity: $\mbX_{s,t} \star \mbX_{t,u}=\mbX_{s,u}$ for all $s\leq t \leq u$,
\item\label{pt:pvar_branched} finite $p$-variation: there exists a control function $\omega$ such that, for all $(s,t)\in\Delta_T$,
\begin{equation}\label{eq:p-var-omega_branched}
\max_{\tau\in F_N}|\scal{\mbX_{s,t},\tau}|^{p/|\tau|} \leq \omega(s,t)\;.
\end{equation}
\end{enumerate}
Properties \ref{pt:Chen_branched}, \ref{pt:pvar_branched}, are obvious analogues of \ref{pt:Chen}, \ref{pt:pvar} for $p$-rough paths.
The property $\mbX_{s,t}\in B^N$ is the substitute for \ref{pt:invert} and the (weak) geometricity condition \eqref{eq:X_shuffle}.

There is a canonical embedding of Hopf algebras $T^N(\R^d)\hookrightarrow H^N_\GL$
given by extending the identity map $\R^d\to\R^d$ to a graded algebra morphism, where we treat 
$e_i\in\R^d$ as an element of
$H^N_\GL$ by $\scal{e_i,\bullet_i}=1$ and $\scal{e_i,\tau}=0$ for all $\tau\in F_N\setminus\{\bullet_i\}$.
This embedding maps $G^N(\R^d)$ into $B^N$ and thus weakly
geometric rough paths embed canonically into branched rough paths.
However, the space of branched rough paths is strictly larger since $B^N$ is larger than the image of $G^N(\R^d)$ under the embedding $T^N(\R^d)\hookrightarrow H^N_\GL$.
To see this in another way, branched rough paths do not impose the shuffle identities \eqref{eq:X_shuffle}.
For example, if $\mbX$ is a geometric (branched) rough path, then
\[
\scal{\mbX,\bullet_i}\scal{\mbX,\bullet_j} =
\scal{\mbX,
\begin{tikzpicture}[scale=0.2,baseline=0.1cm]
        \node at (0,0)  [dot,label= {[label distance=-0.2em]right: \scriptsize  $ i $} ] (root) {};
         \node at (0,2)  [dot,label={[label distance=-0.2em]right: \scriptsize  $ j$}] (right) {};
            \draw[kernel1] (right) to
     node [sloped,below] {\small }     (root);
     \end{tikzpicture}
}
+
\scal{\mbX,
\begin{tikzpicture}[scale=0.2,baseline=0.1cm]
        \node at (0,0)  [dot,label= {[label distance=-0.2em]right: \scriptsize  $ j $} ] (root) {};
         \node at (0,2)  [dot,label={[label distance=-0.2em]right: \scriptsize  $ i$}] (right) {};
            \draw[kernel1] (right) to
     node [sloped,below] {\small }     (root);
     \end{tikzpicture}
}\;,
\]
while this need not hold for general branched rough paths.

To see how non-geometric branched rough paths may arise, consider $X\colon[0,T]\to\R^d$ a standard $d$-dimensional Brownian motion
and define $\mbX$ as in \eqref{eq:II_branched}  where the integral is taken in the It\^o sense.
Then, by It\^o's formula,
\[
\scal{\mbX_{0,t},\bullet_i}\scal{\mbX_{0,t},\bullet_j} =
\int_0^t X_s^i\mrd X^j_s + \int_0^t X_s^j\mrd X^i_s + \delta_{ij} t/2\;,
\]
where the first two integrals on the right-hand side are
$
\scal{\mbX_{0,t},
\begin{tikzpicture}[scale=0.2,baseline=0.1cm]
        \node at (0,0)  [dot,label= {[label distance=-0.2em]right: \scriptsize  $ i $} ] (root) {};
         \node at (0,2)  [dot,label={[label distance=-0.2em]right: \scriptsize  $ j$}] (right) {};
            \draw[kernel1] (right) to
     node [sloped,below] {\small }     (root);
     \end{tikzpicture}
}
$
and
$
\scal{\mbX_{0,t},
\begin{tikzpicture}[scale=0.2,baseline=0.1cm]
        \node at (0,0)  [dot,label= {[label distance=-0.2em]right: \scriptsize  $ j $} ] (root) {};
         \node at (0,2)  [dot,label={[label distance=-0.2em]right: \scriptsize  $ i$}] (right) {};
            \draw[kernel1] (right) to
     node [sloped,below] {\small }     (root);
     \end{tikzpicture}
}
$
respectively.
Therefore $\mbX$, which is almost surely a $1/p$-H\"older branched rough path for any $p>2$ (i.e. we can take $\omega(s,t)=|t-s|$ in \eqref{eq:p-var-omega_branched})
is not a geometric rough path.
(If we chose Stratonovich instead of It\^o integration, the final term $\delta_{ij}t/2$ would not arise and we would have obtained a geometric rough path as before.)
This additional flexibility does not only allow for a more canonical way to treat It\^o integrals in a rough path sense (without needing to convert between It\^o and Stratonovich),
but is a crucial ingredient in linking rough paths and regularity structures that we discuss below,
since renormalisation procedures in the latter break the classical chain rule that the shuffle identities imply.

The notion of a controlled rough path extends almost verbatim to the branched setting and allows for a generalisation of the ULT.
A comparison between the two approaches, including a characterisation of solutions to RDEs based on Euler estimates akin to Davie are discussed by \cite{HairerKelly15}.
The authors therein also devise a (non-canonical) way to move from non-geometric to geometric rough paths by means of the Lyons--Victoir extension theorem.

An extension of a branched rough path $\mbX$ is a map $\hat\mbX\colon \Delta_T \to H^M_\GL$ for $M > N$
such that $\scal{\hat\mbX,\tau}=\scal{\mbX,\tau}$ for all $|\tau|\leq N$
and which satisfies \eqref{eq:p-var-omega_branched} with $N$ replaced by $M$.
Extensions to any level $M> N$ for branched rough paths exist as in the geometric case, which can be shown inductively by integrating $\mbX$ against itself as a controlled path.
In contrast to the geometric setting, however, one cannot in general view $\hat\mbX$ as the solution to a linear RDE driven by $\mbX$ as in \eqref{eq:RDE_extension}; indeed, there exists non-constant (and necessarily non-geometric) branched rough paths $\mbX$ such that every linear RDE driven by $\mbX$ is constant.
Analogues of the decay estimates \eqref{eq:p_var_sharp} for branched rough path extensions were shown by \cite{Boedihardjo18}.

Note that $T^2(\R^d)$ is canonically isomorphic to $H^2_\GL$,
which means that (a version of) the ULT holds for $p$-rough paths in the sense of Definition \ref{def:RP} for $p\in [2,3)$.
Lejay--Victoir found a canonical isomorphism for $p\in (2,3)$ between (branched) $p$-rough paths and geometric rough paths with mixed $(p,p/2)$-variation over a larger vector space.
By using an isomorphism of Foissy--Chapoton between $H_\GL$ and a tensor algebra over an infinite dimensional vector space,
this canonical isomorphism between branched $p$-rough paths and geometric rough paths was extended by \cite{BC19} to all $p\geq 1$, which in particular extends the results in \citep{BGLY16,CL16} to branched rough paths.

\textbf{Discontinuous processes.}
There exist extensions of rough path theory to the case where $\mbX$ is discontinuous.
In the bounded variation case, there are (at least) two distinct ways to interpret the ODE \eqref{eq:ODE} when $X$ is discontinuous.
The first is through an integral equation $Y_t=Y_0+\int_0^t Y_s\mrd X_s$, where the integral is defined as the limit of Riemann sums as in \eqref{eq:Y_X_int};
we call this a forward (or It{\^o}) equation/integral.
The second is by interpreting the jumps of $X$ as infinitesimal continuous motions, which, unlike the forward case, respects the underlying geometry of the vector fields $f$ in \eqref{eq:ODE} (if $f$ takes values in the tangent space of a submanifold of $W$, then the solution $Y$ in the geometric sense will remain inside the submanifold).
These two approaches both have meaningful generalisations to the case of rough paths.

\begin{figure}[t]
\centering
\begin{subfigure}[h]{0.25\textwidth}
\begin{tikzpicture}[scale=0.8]
          \tikzmath{ \hhh = 3.14/2; \hhhh = \hhh * 1.2; }
          \begin{axis}[
unit vector ratio*=1 1 1,
			  disabledatascaling,
              height       = 1.8in,
              xmax         = 1.7,
              ymax         = \hhhh ,
ymin         = -0.1 ,
              xtick        = {0.0, 1},
              xticklabels  = {0, $1$ },
              ytick        = {0.0, \hhh },
              yticklabels  = {0, $\pi/2$},
              axis lines   = center,
              line cap=round,
x label style={at={(axis description cs:0.9,-0.0)},anchor=north},
    y label style={at={(axis description cs:-0.1,0.3)},anchor=south},
    axis line style={->},
xlabel= {$t$},
ylabel= {$X_{t}$}
              ]
\node at (1,\hhh) [circle,color=red,fill,inner sep=1.5pt]{};
              \addplot[color=gray,dotted] coordinates { (0,\hhh) (1,\hhh) };
              \addplot[color=red,very thick] coordinates { (0,0) (1,0) };
              \addplot[color=red,very thick] coordinates { (1, \hhh ) (2, \hhh ) };
          \end{axis}
      \end{tikzpicture}
\caption{}
\end{subfigure}
\begin{subfigure}[h]{0.35\textwidth}
      \begin{tikzpicture}[scale=0.8]
          \tikzmath{ \hhh = 3.14/2;  \hhhh = \hhh * 1.5;}
          \begin{axis}[
unit vector ratio*=1 1 1,
              line cap   = round,
              view       = {-45}{-25},
              axis lines = center,
              xmax       = 1.5,
              ymax       = \hhhh,
              zmax       = 2.1,
              height     = 3in,
              xtick      = {1},
              ytick      = {\hhh},
              yticklabel = {$\pi/2$},
              ztick      = {1},
              zticklabel = {1},
x label style={at={(axis description cs:1,0.4)},anchor=north},
xlabel= {$Y^1$},
ylabel= {$Y^2$},
zlabel={$t$}
              ]
\addplot3[mark=*,red,point meta=explicit symbolic,nodes near coords] coordinates {(1,\hhh,1)};
              \addplot3[
                  color      = red,
                  mark       = none,
                  very thick
                  ] coordinates { (1,0,0) (1,0,1)};
\addplot3[dotted,
                  color      = gray,
                  mark       = none,
                   thick
                  ] coordinates { (1,0,1) (0,0,1)};
\addplot3[dotted,
                  color      = gray,
                  mark       = none,
                   thick
                  ] coordinates { (1,\hhh,0) (1,0,0)};
\addplot3[dotted,
                  color      = gray,
                  mark       = none,
                   thick
                  ] coordinates { (1,\hhh,0) (0,\hhh,0)};
\addplot3[dotted,
                  color      = gray,
                  mark       = none,
                   thick
                  ] coordinates { (1,\hhh,1) (1,\hhh,0)};
\addplot3[dotted,
                  color      = gray,
                  mark       = none,
                   thick
                  ] coordinates { (1,\hhh,1) (1,0,1)};
              \addplot3[
                  color      = red,
                  mark       = none,
                  very thick
                  ] coordinates { (1,\hhh,1)
                  (1,\hhh,2.0) };
          \end{axis}
      \end{tikzpicture}
\caption{}
\end{subfigure}
\begin{subfigure}[h]{0.35\textwidth}
\begin{tikzpicture}[scale=0.8]
          \tikzmath{ \hhh = 3.14/2;}
          \begin{axis}[
unit vector ratio*=1 1 1,
              line cap   = round,
              view       = {-45}{-25},
              axis lines = center,
              xmax       = 1.5,
              ymax       = 1.5,
              zmax       = 1.7,
              height     = 3in,
              xtick      = {1},
              ytick      = {1},
              yticklabel = {$1$},
              ztick      = {1},
              zticklabel = {1},
x label style={at={(axis description cs:1,0.4)},anchor=north},
xlabel= {$Y^1$},
ylabel= {$Y^2$},
zlabel={$t$}
              ]
              \addplot3+ [
                  domain     = 0:\hhh ,
                  samples    = 64,
                  samples y  = 0,
                  mark       = none,
                  very thick,
                  color      = red
                  ] ( {cos(deg(x))},{sin(deg(x))},{1 + 0.0 * x / \hhh } );
              \addplot3+ [
                  color      = red,
                  mark       = none,
                  very thick
                  ] coordinates { (1,0,0) (1,0,1)};
              \addplot3+ [
                  color      = red,
                  mark       = none,
                  very thick
                  ] coordinates { ({cos(deg(\hhh))},{sin(deg(\hhh))},1)
                  ({cos(deg(\hhh))},{sin(deg(\hhh))},2) };
\addplot3[dotted,
                  color      = gray,
                  mark       = none,
                   thick
                  ] coordinates { (1,0,1) (0,0,1)};
\addplot3[dotted,
                  color      = gray,
                  mark       = none,
                   thick
                  ] coordinates { (0,1,1) (0,0,1)};
\addplot3[dotted,
                  color      = gray,
                  mark       = none,
                   thick
                  ] coordinates { (0,1,1) (0,1,0)};
          \end{axis}
      \end{tikzpicture}
\caption{}
\end{subfigure}
\caption{(a): the path $X_t = \frac\pi2\bone_{t\geq 1}$. (b) and (c): solutions to $\begin{pmatrix}
      \mrd Y^1 _t\\ \mrd Y^2_t
    \end{pmatrix}
    =
    \begin{pmatrix}
      - Y_t^2 \\ Y^1_t
    \end{pmatrix}
    \mrd X_t$ in forward and geometric sense respectively.}
\label{fig:compare}
\end{figure}
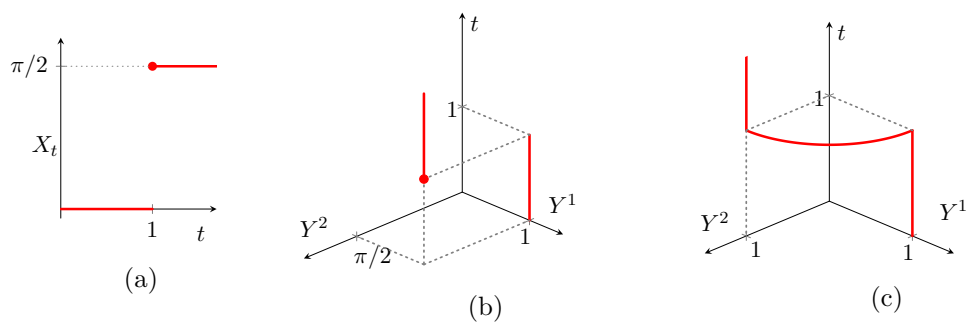

The geometric approach, inspired by earlier work of Marcus, was first extended to the rough case by \cite{Williams01} who restricted attention to L\'evy processes.
To describe the idea, suppose $X\colon[0,T]\to \R^d$ has bounded variation and is c\`adl\`ag (continue \`a droite, limite \`a gauche, i.e. right continuous with left limits)
and that we wish to give meaning to a solution of $\mrd Y = f(Y)\mrd X$.
At every jump time $t$ of $X$,
we add an `infinitesimal time' $\delta t$ over which $X$ traverses its jump by a continuous path $X^{(t)}$.
The solution $Y$ correspondingly flows along the vector fields $f$ as if it is driven by $X^{(t)}$ in time $\delta t$.
In the approach of Marcus and Williams, $X^{(t)}$ is taken linear, see Figure \ref{fig:Marcus}.
It is important to note that the flow of $Y$ over $\delta t$ need not be linear
and that this notion of solution is different from interpreting the ODE as an integral equation $Y_t = Y_0+\int_0^t f(Y_s)\mrd X_s$, see Figure \ref{fig:compare}.

\begin{figure}[t]
\centering
\begin{tikzpicture}[scale=0.7]

\begin{axis}[
  xtick={0,4,8}, ytick={0,4,8},  xticklabels={,,},  yticklabels={,,}, xlabel style={below right}, ylabel style={above left},
  xmin= -3, xmax=10,  ymin=-3,  ymax=10]

\addplot[domain=0:3,blue, thick] ({-1+2*x},{7-(x-3)^2});   
\addplot[domain=3:6,blue, thick, samples = 200] ({ 7- 2*(x-3)},{4 + 1.5*sin(  (x-3)/2 * 3 *360 ) }); 
\addplot[domain=6:12,blue, thick, samples= 200] ({0 - (x-6)/2},{8-(x-6)/4 - 2*sin(  (x-6)/(10-6) * 1 *360 ) / 3 });
\draw[dotted, thick, red,-{Latex[length=3mm]}] (axis cs:5,7) -- (axis cs:7,4);
\draw[dotted, thick, red,-{Latex[length=3mm]}] (axis cs:1,4) -- (axis cs:0,8);

\addplot[soldot] coordinates{(7,4)(0,8)};   
\addplot[holdot] coordinates{(5,7)(1,4)};   

\end{axis}
\end{tikzpicture}
\qquad
\begin{tikzpicture}[scale=0.7]

\begin{axis}[
  xtick={0,4,8}, ytick={0,4,8},  xticklabels={,,},  yticklabels={,,}, xlabel style={below right}, ylabel style={above left},
  xmin= -3, xmax=10,  ymin=-3,  ymax=10]

\addplot[domain=0:3,blue, thick, samples= 200] ({-1+2*x},{7-(x-3)^2+sin(x *360) +0.5*(3-x)});
\addplot[domain=3:6,blue, thick, samples = 200] ({ 7- 2*(x-3)},{4 + sin(  (x-3)/2 * 3 *360 ) -2*(x-3)^2/3}); 
\addplot[domain=6:12,blue, thick, samples= 200] ({0 - (x-6)/2},{8+(x-6)/4 - 2*sin(  (x-6)/(10-6) * 1 *360 ) / 2 });

\addplot[domain=0:1,red,dotted, thick, samples= 200,-{Latex[length=3mm]}] ({5+2*x - sin(x*360)},{7-3*x});

\addplot[domain=0:1,red,dotted, thick, samples= 200,-{Latex[length=3mm]}] ({1-x+sin(x*360)},{-2+10*x});

\addplot[soldot] coordinates{(7,4)(0,8)};   
\addplot[holdot] coordinates{(5,7)(1,-2)};   

\end{axis}
\end{tikzpicture}
\caption{Linear (Marcus) interpolation of c\`adl\`ag driver $X$ (left) and the corresponding geometric solution to $\mrd Y = f(Y)\mrd X$ (right).}
\label{fig:Marcus}
\end{figure}
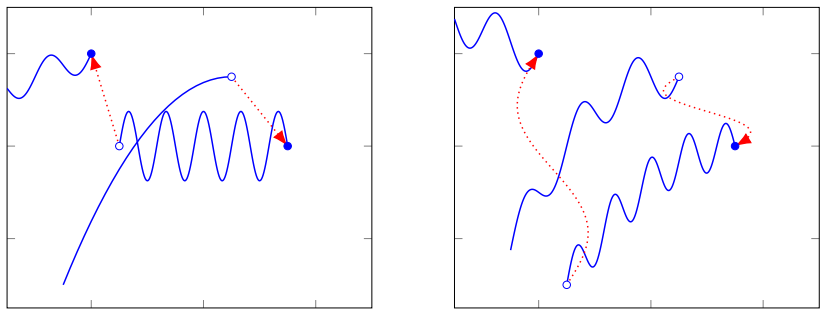

This geometric approach was developed further in the rough path setting by several authors.
\cite{FS17} initiated the systematic study of rough integration and RDEs with jumps (both in geometric and forward sense).
\cite{Chevyrev18} characterised all possible $G^N(\R^d)$-valued L\'evy processes that give rise to $p$-rough paths.
\cite{CF19} extended the ULT to geometric rough paths with jumps;
a feature of the metric appearing in this ULT is that
it requires working with so-called path functions from \citep{Chevyrev18}
that allow for more general paths to connect the jumps of $X$ than straight lines.
Furthermore smooth paths are dense in the space of c\`adl\`ag rough paths under this metric, which is not the case for the classical Skorokhod $J_1$ topology.
\cite{CF19} also showed a rough path version of L\'epingle's $p$-variation Burkholder--Davis--Gundy (BDG) inequality, $\E \|\mbX\|_{\var p}^{q} \asymp \E[X]^{q/2}$ for $p>2$ and $q\geq 1$, where $[X]$ is the quadratic variation of a semi-martingale $X$ with canonical lift $\mbX$.
This BDG inequality and ULT imply stability results for Marcus SDEs under the uniformly-controlled-variations condition of Kurtz--Protter and Jakubowski--M\'emin--Pag\`es (cf. \cite{Kurtz_Protter_91} for such a result for It\^o SDEs)
as well as extend the Wong--Zakai theorem of \cite{KPP95}.

The ULT for forward/It\^o rough differential equations was shown by \cite{FZ18}.
Unlike the geometric case where the proof of the ULT follows from the continuous theory,
the forward case requires reworking many of the analytic steps
(including a sewing lemma for discontinuous control functions).
In addition, due to the non-geometric nature of the resulting RDEs (in particular the absence of the classical chain rule), the natural state space to work with for $p\geq 2$ is the GL algebra as in the case of branched rough paths.

A problem in which both forward and geometric solution theories arise is homogenisation of deterministic fast-slow systems:
the forward approach was used by \cite{CFKMZ19} to extend the results of \cite{KM17} to a discrete-time setting,
while the geometric approach was used by \cite{CFKM20} to prove homogenisation for superdiffusive fast-slow systems.

\textbf{Singular stochastic partial differential equations.}
Over the past decade there has been much progress in the analysis of stochastic partial differential equations (SPDEs)
using ideas coming from rough path theory.
Notably, rough paths in the form given above (albeit with non-trivial renormalisation procedures)
were used by \cite{Hairer_13_KPZ} to solve the Khardar--Parisi--Zhang (KPZ) equation
\begin{equation}\label{eq:KPZ}
\partial_t u = \partial_x^2 u + (\partial_x u)^2 + \xi 
\end{equation}
where $u\colon [0,T]\times \R\to\R$ and $\xi$ is a space-time white noise on $\R\times \R$.
To see why a solution theory for \eqref{eq:KPZ} is non-trivial,
we note that an argument of Kolmogorov implies that $\xi$ is almost surely in the H\"older--Besov space of negative regularity $C^{\alpha}$ for $\alpha<-\frac32$ measured in the parabolic scaling, but not for $\alpha\geq -\frac32$.
Consequently, we expect (and can ultimately show) that $u$ has H\"older regularity at best $\alpha+2<\frac12$, where the $+2$ comes from the heat operator $\partial_t - \partial_x^2$.
However, this means that $\partial_x u$ is a distribution of regularity $\alpha-1<-\frac12$, which renders its square $(\partial_x u)^2$ analytically ill-defined.
For this final point,
we recall that the product map $(f,g)\mapsto fg$, defined for smooth $f,g$, extends to a continuous bilinear map $C^\alpha\times C^\beta \to C^\alpha$ for $\alpha\leq \beta$ if and only if $\alpha+\beta>0$.
This shows there is little hope to naively solve \eqref{eq:KPZ} by using, e.g. Picard iterations.

This type of powering counting is of course applicable to classical SDEs $\partial_t Y = f(Y)\xi$ where $\xi = \mrd X_t/\mrd t$ is a white noise (time derivative of a Brownian motion $X$):
we know $\xi$ is at best in $C^\alpha$ for $\alpha<-\frac12$, so $u$ and thus $f(u)$ are at best in $C^{\alpha+1}$, but this renders the product $f(u)\xi$ analytically ill-posed since $(\alpha+1)+\alpha<0$.
This partially justifies why stochastic calculus falls outside the scope of deterministic ODE theory (\cite{Lyons91} justifies this in a stronger sense as discussed above).
Nonetheless, one can give meaning to $\partial_t Y = f(Y)\xi$ by taking a mollified approximation $\xi^\eps$ to $\xi$ and showing that the corresponding solutions $Y^\eps$ converge to a limit $Y$
- this is the Wong--Zakai theorem that furthermore identifies $Y$ with the solution of the SDE in the Stratonovich sense.

In contrast, there is a further complication with \eqref{eq:KPZ} in that, if we take mollifications $\xi^\eps$ of $\xi$
and denote by $u^\eps$ the  corresponding solutions to \eqref{eq:KPZ},
then $u^\eps$ \textit{diverges} as $\eps\to 0$.
Instead, it turns out one needs to insert a counterterm $-C_\eps\in \R$ to the right-hand side of \eqref{eq:KPZ} with $C_\eps\to\infty$ as $\eps\to0$ in order for $u^\eps$ to converge to a limit, as shown by \cite{Hairer_13_KPZ}.

Shortly after, \cite{Hairer14} devised the theory of regularity structures
that is a generalisation of rough paths to higher dimensions and is able to give a robust solution theory for a wide class of SPDEs.
Prototypical equations that fall within its scope are non-linear heat equations posed on $[0,T]\times \R^n$ of the form
\begin{equation}\label{eq:SPDE}
\partial_t u = \Delta u + f(u,\nabla u ,\xi)\;,
\end{equation}
where $\xi$ is a distribution and $f$ is a smooth function affine in $\xi$.
Several technical assumptions are needed on $f$ and $\xi$, the most fundamental one of which is sub-criticality,
which roughly states that $u$ is a perturbation of the solution of the linearised equation $\partial_t v = \Delta v+\xi$.
In addition to \eqref{eq:KPZ}, examples that falls within the scope of regularity structures are the parabolic Anderson model in spatial dimensions $2$ and $3$ and the parabolic stochastic quantisation equation of $\Phi^4_3$ model in quantum field theory.

The philosophy in regularity structures is close to that of controlled rough paths:
one enhances the original signal $\xi$ to a richer object called a \textit{model} (analogue of a rough path)
and solves a `lifted' equation in a space of \textit{modelled distributions} (analogue of controlled path).
As in rough path theory, this solution map model $\mapsto$ modelled distribution is continuous.
%

In parallel to the development of regularity structures, \cite{GIP15} introduced the theory of paracontrolled distributions, the purpose of which is also to give a 
pathwise solution theory for singular SPDEs of the type \eqref{eq:SPDE}.
The philosophy in this approach is again to lift $\xi$ to an enhanced data set followed by a continuous solution map on functions satisfying a `paracontrolled ansatz' (the analogue of a controlled path/modelled distribution).
In contrast to regularity structures that are based on local expansions of jets, paracontrolled distributions are Fourier analytic relying on Bony's paraproduct.

A typical output of either theory is that solutions to \eqref{eq:SPDE}, with $\xi$ replaced by a mollification $\xi^\eps$, converge provided that a finite number of suitable counterterms ($\eps$-dependent functions of $u$) are added to the right-hand side of \eqref{eq:SPDE}.
These counterterms arise from the need to `renormalise' the canonical enhancement of $\xi^\eps$,
which then has a corresponding effect on the equation; these renormalisation procedures are carried out by \cite{BHZ16,BCCH21} in the framework of regularity structures
with the use of Hopf and pre-Lie algebras that generalise the Connes--Kreimer algebra.
Renormalisation, although one of the main differences in how regularity structures and rough paths are applied in practice,
can also be understood as translations of branched rough paths~\citep[see][]{BCFP19}.

\bibliography{./refs}

\def\cprime{$'$} \def\polhk#1{\setbox0=\hbox{#1}{\ooalign{\hidewidth
  \lower1.5ex\hbox{`}\hidewidth\crcr\unhbox0}}}
\begin{thebibliography}{}

\bibitem[Bass et~al., 2002]{BHL02}
Bass, R.~F., Hambly, B.~M., and Lyons, T.~J. (2002).
\newblock Extending the {W}ong-{Z}akai theorem to reversible {M}arkov
  processes.
\newblock {\em J. Eur. Math. Soc. (JEMS)}, 4(3):237--269.

\bibitem[Bestvina, 2002]{Bestvina_02_trees}
Bestvina, M. (2002).
\newblock {$\Bbb R$}-trees in topology, geometry, and group theory.
\newblock In {\em Handbook of geometric topology}, pages 55--91. North-Holland,
  Amsterdam.

\bibitem[Boedihardjo, 2018]{Boedihardjo18}
Boedihardjo, H. (2018).
\newblock Decay rate of iterated integrals of branched rough paths.
\newblock {\em Ann. Inst. H. Poincar\'{e} C Anal. Non Lin\'{e}aire},
  35(4):945--969.

\bibitem[Boedihardjo and Chevyrev, 2019]{BC19}
Boedihardjo, H. and Chevyrev, I. (2019).
\newblock An isomorphism between branched and geometric rough paths.
\newblock {\em Ann. Inst. Henri Poincar\'{e} Probab. Stat.}, 55(2):1131--1148.

\bibitem[Boedihardjo et~al., 2016]{BGLY16}
Boedihardjo, H., Geng, X., Lyons, T., and Yang, D. (2016).
\newblock The signature of a rough path: uniqueness.
\newblock {\em Adv. Math.}, 293:720--737.

\bibitem[Bonetto et~al., 2019]{BCKL19}
Bonetto, F., Chernov, N., Korepanov, A., and Lebowitz, J.~L. (2019).
\newblock Autonomous evolution of electron speeds in a thermostatted system:
  exact results.
\newblock {\em Nonlinearity}, 32(6):2055--2072.

\bibitem[Broux and Zambotti, 2022]{BZ22}
Broux, L. and Zambotti, L. (2022).
\newblock The sewing lemma for {$0<\gamma\leq1$}.
\newblock {\em J. Funct. Anal.}, 283(10):Paper No. 109644, 34.

\bibitem[Bruned et~al., 2021]{BCCH21}
Bruned, Y., Chandra, A., Chevyrev, I., and Hairer, M. (2021).
\newblock Renormalising {SPDE}s in regularity structures.
\newblock {\em J. Eur. Math. Soc. (JEMS)}, 23(3):869--947.

\bibitem[Bruned et~al., 2019a]{BCFP19}
Bruned, Y., Chevyrev, I., Friz, P.~K., and Prei\ss, R. (2019a).
\newblock A rough path perspective on renormalization.
\newblock {\em J. Funct. Anal.}, 277(11):108283, 60.

\bibitem[Bruned et~al., 2019b]{BHZ16}
Bruned, Y., Hairer, M., and Zambotti, L. (2019b).
\newblock Algebraic renormalisation of regularity structures.
\newblock {\em Invent. Math.}, 215(3):1039--1156.

\bibitem[Cass et~al., 2016]{Cass_Driver_Lim_16_wgrp}
Cass, T., Driver, B.~K., Lim, N., and Litterer, C. (2016).
\newblock On the integration of weakly geometric rough paths.
\newblock {\em J. Math. Soc. Japan}, 68(4):1505--1524.

\bibitem[Cass and Friz, 2010]{CF10}
Cass, T. and Friz, P. (2010).
\newblock Densities for rough differential equations under {H}\"{o}rmander's
  condition.
\newblock {\em Ann. of Math. (2)}, 171(3):2115--2141.

\bibitem[Cass et~al., 2015]{CHLT15}
Cass, T., Hairer, M., Litterer, C., and Tindel, S. (2015).
\newblock Smoothness of the density for solutions to {G}aussian rough
  differential equations.
\newblock {\em Ann. Probab.}, 43(1):188--239.

\bibitem[Cass and Lyons, 2015]{Cass_Lyons_15}
Cass, T. and Lyons, T. (2015).
\newblock Evolving communities with individual preferences.
\newblock {\em Proc. Lond. Math. Soc. (3)}, 110(1):83--107.

\bibitem[Chevyrev, 2018]{Chevyrev18}
Chevyrev, I. (2018).
\newblock Random walks and {L}\'{e}vy processes as rough paths.
\newblock {\em Probab. Theory Related Fields}, 170(3-4):891--932.

\bibitem[Chevyrev et~al., 2022]{CFKMZ19}
Chevyrev, I., Friz, P., Korepanov, A., Melbourne, I., and Zhang, H. (2022).
\newblock Deterministic homogenization under optimal moment assumptions for
  fast--slow systems. {P}art 2.
\newblock {\em Ann. Inst. Henri Poincar\'{e} Probab. Stat.}, 58(3):1328--1350.

\bibitem[Chevyrev and Friz, 2019]{CF19}
Chevyrev, I. and Friz, P.~K. (2019).
\newblock Canonical {RDE}s and general semimartingales as rough paths.
\newblock {\em Ann. Probab.}, 47(1):420--463.

\bibitem[Chevyrev et~al., 2020]{CFKM20}
Chevyrev, I., Friz, P.~K., Korepanov, A., and Melbourne, I. (2020).
\newblock Superdiffusive limits for deterministic fast-slow dynamical systems.
\newblock {\em Probab. Theory Related Fields}, 178(3-4):735--770.

\bibitem[{Chevyrev} and {Kormilitzin}, 2016]{CK16}
{Chevyrev}, I. and {Kormilitzin}, A. (2016).
\newblock {A Primer on the Signature Method in Machine Learning}.
\newblock {\em arXiv e-prints}.

\bibitem[Chevyrev and Lyons, 2016]{CL16}
Chevyrev, I. and Lyons, T. (2016).
\newblock Characteristic functions of measures on geometric rough paths.
\newblock {\em Ann. Probab.}, 44(6):4049--4082.

\bibitem[Chevyrev and Oberhauser, 2022]{CO22}
Chevyrev, I. and Oberhauser, H. (2022).
\newblock Signature moments to characterize laws of stochastic processes.
\newblock {\em Journal of Machine Learning Research}, 23(176):1--42.

\bibitem[Coghi et~al., 2020]{CDFM20}
Coghi, M., Deuschel, J.-D., Friz, P.~K., and Maurelli, M. (2020).
\newblock Pathwise {M}c{K}ean-{V}lasov theory with additive noise.
\newblock {\em Ann. Appl. Probab.}, 30(5):2355--2392.

\bibitem[Coutin and Qian, 2002]{Coutin_Qian_02}
Coutin, L. and Qian, Z. (2002).
\newblock Stochastic analysis, rough path analysis and fractional {B}rownian
  motions.
\newblock {\em Probab. Theory Related Fields}, 122(1):108--140.

\bibitem[Davie, 2008]{Davie08}
Davie, A.~M. (2008).
\newblock Differential equations driven by rough paths: an approach via
  discrete approximation.
\newblock {\em Appl. Math. Res. Express. AMRX}, 2008:Art. ID abm009, 40.

\bibitem[Deya et~al., 2019]{DGHT19}
Deya, A., Gubinelli, M., Hofmanov\'{a}, M., and Tindel, S. (2019).
\newblock A priori estimates for rough {PDE}s with application to rough
  conservation laws.
\newblock {\em J. Funct. Anal.}, 276(12):3577--3645.

\bibitem[{Fermanian} et~al., 2023]{FLMS_23}
{Fermanian}, A., {Lyons}, T., {Morrill}, J., and {Salvi}, C. (2023).
\newblock {New directions in the applications of rough path theory}.
\newblock {\em arXiv e-prints}.
\newblock To appear in \textit{IEEE BITS the Information Theory Magazine}.

\bibitem[Feyel and de~La~Pradelle, 2006]{Feyel_de_La_Pradelle_06}
Feyel, D. and de~La~Pradelle, A. (2006).
\newblock Curvilinear integrals along enriched paths.
\newblock {\em Electron. J. Probab.}, 11:no. 34, 860--892.

\bibitem[Friz and Victoir, 2008]{FV08}
Friz, P. and Victoir, N. (2008).
\newblock On uniformly subelliptic operators and stochastic area.
\newblock {\em Probab. Theory Related Fields}, 142(3-4):475--523.

\bibitem[Friz and Hairer, 2020]{FrizHairer20}
Friz, P.~K. and Hairer, M. (2020).
\newblock {\em A course on rough paths}.
\newblock Universitext. Springer, Cham.
\newblock With an introduction to regularity structures, Second edition of [
  3289027].

\bibitem[Friz and Shekhar, 2017]{FS17}
Friz, P.~K. and Shekhar, A. (2017).
\newblock General rough integration, {L}\'evy rough paths and a
  {L}\'evy-{K}intchine-type formula.
\newblock {\em Ann. Probab.}, 45(4):2707--2765.

\bibitem[Friz and Victoir, 2010]{FV10}
Friz, P.~K. and Victoir, N.~B. (2010).
\newblock {\em Multidimensional Stochastic Processes as Rough Paths: Theory and
  Applications}.
\newblock Cambridge Studies in Advanced Mathematics. Cambridge University
  Press.

\bibitem[Friz and Zhang, 2018]{FZ18}
Friz, P.~K. and Zhang, H. (2018).
\newblock Differential equations driven by rough paths with jumps.
\newblock {\em J. Differential Equations}, 264(10):6226--6301.

\bibitem[Grong et~al., 2022]{Grong_Nilssen_Schmedding_22_wgrp}
Grong, E., Nilssen, T., and Schmeding, A. (2022).
\newblock Geometric rough paths on infinite dimensional spaces.
\newblock {\em J. Differential Equations}, 340:151--178.

\bibitem[Gubinelli, 2004]{Gubinelli04}
Gubinelli, M. (2004).
\newblock Controlling rough paths.
\newblock {\em J. Funct. Anal.}, 216(1):86--140.

\bibitem[Gubinelli, 2010]{Gubinelli_10}
Gubinelli, M. (2010).
\newblock Ramification of rough paths.
\newblock {\em J. Differential Equations}, 248(4):693--721.

\bibitem[Gubinelli et~al., 2015]{GIP15}
Gubinelli, M., Imkeller, P., and Perkowski, N. (2015).
\newblock Paracontrolled distributions and singular {PDE}s.
\newblock {\em Forum Math. Pi}, 3:e6, 75.

\bibitem[Hairer, 2013]{Hairer_13_KPZ}
Hairer, M. (2013).
\newblock Solving the {KPZ} equation.
\newblock {\em Ann. of Math. (2)}, 178(2):559--664.

\bibitem[Hairer, 2014]{Hairer14}
Hairer, M. (2014).
\newblock A theory of regularity structures.
\newblock {\em Invent. Math.}, 198(2):269--504.

\bibitem[Hairer and Kelly, 2015]{HairerKelly15}
Hairer, M. and Kelly, D. (2015).
\newblock Geometric versus non-geometric rough paths.
\newblock {\em Ann. Inst. Henri Poincar\'e Probab. Stat.}, 51(1):207--251.

\bibitem[Hambly and Lyons, 2010]{HL10}
Hambly, B. and Lyons, T. (2010).
\newblock Uniqueness for the signature of a path of bounded variation and the
  reduced path group.
\newblock {\em Ann. of Math. (2)}, 171(1):109--167.

\bibitem[Hara and Hino, 2010]{Hara_Hino_10}
Hara, K. and Hino, M. (2010).
\newblock Fractional order {T}aylor's series and the neo-classical inequality.
\newblock {\em Bull. Lond. Math. Soc.}, 42(3):467--477.

\bibitem[Ikeda and Watanabe, 1989]{Ikeda_Watanabe_89_SDEs}
Ikeda, N. and Watanabe, S. (1989).
\newblock {\em Stochastic differential equations and diffusion processes},
  volume~24 of {\em North-Holland Mathematical Library}.
\newblock North-Holland Publishing Co., Amsterdam; Kodansha, Ltd., Tokyo,
  second edition.

\bibitem[Kelly and Melbourne, 2016]{KM16}
Kelly, D. and Melbourne, I. (2016).
\newblock Smooth approximation of stochastic differential equations.
\newblock {\em Ann. Probab.}, 44(1):479--520.

\bibitem[Kelly and Melbourne, 2017]{KM17}
Kelly, D. and Melbourne, I. (2017).
\newblock Deterministic homogenization for fast-slow systems with chaotic
  noise.
\newblock {\em J. Funct. Anal.}, 272(10):4063--4102.

\bibitem[Kiraly and Oberhauser, 2019]{KO19}
Kiraly, F.~J. and Oberhauser, H. (2019).
\newblock Kernels for sequentially ordered data.
\newblock {\em Journal of Machine Learning Research}, 20(31):1--45.

\bibitem[Kurtz et~al., 1995]{KPP95}
Kurtz, T.~G., Pardoux, {\'E}., and Protter, P. (1995).
\newblock Stratonovich stochastic differential equations driven by general
  semimartingales.
\newblock {\em Ann. Inst. H. Poincar\'e Probab. Statist.}, 31(2):351--377.

\bibitem[Kurtz and Protter, 1991]{Kurtz_Protter_91}
Kurtz, T.~G. and Protter, P. (1991).
\newblock Weak limit theorems for stochastic integrals and stochastic
  differential equations.
\newblock {\em Ann. Probab.}, 19(3):1035--1070.

\bibitem[Li and Ni, 2022]{Li_Ni_22}
Li, S. and Ni, H. (2022).
\newblock Expected signature of stopped {B}rownian motion on {$d$}-dimensional
  {$C^{2,\alpha}$}-domains has finite radius of convergence everywhere: {$2\leq
  d\leq8$}.
\newblock {\em J. Funct. Anal.}, 282(12):Paper No. 109447, 45.

\bibitem[Lyons, 1991]{Lyons91}
Lyons, T. (1991).
\newblock On the nonexistence of path integrals.
\newblock {\em Proc. Roy. Soc. London Ser. A}, 432(1885):281--290.

\bibitem[Lyons and Qian, 2002]{LyonsQian02}
Lyons, T. and Qian, Z. (2002).
\newblock {\em System control and rough paths}.
\newblock Oxford Mathematical Monographs. Oxford University Press, Oxford.
\newblock Oxford Science Publications.

\bibitem[Lyons and Victoir, 2007]{Lyons_Victoir_07}
Lyons, T. and Victoir, N. (2007).
\newblock An extension theorem to rough paths.
\newblock {\em Ann. Inst. H. Poincar\'e Anal. Non Lin\'eaire}, 24(5):835--847.

\bibitem[Lyons, 1998]{Lyons98}
Lyons, T.~J. (1998).
\newblock Differential equations driven by rough signals.
\newblock {\em Rev. Mat. Iberoamericana}, 14(2):215--310.

\bibitem[Lyons and Xu, 2018]{LyonsXu18}
Lyons, T.~J. and Xu, W. (2018).
\newblock Inverting the signature of a path.
\newblock {\em J. Eur. Math. Soc. (JEMS)}, 20(7):1655--1687.

\bibitem[Williams, 2001]{Williams01}
Williams, D. R.~E. (2001).
\newblock Path-wise solutions of stochastic differential equations driven by
  {L}\'evy processes.
\newblock {\em Rev. Mat. Iberoamericana}, 17(2):295--329.

\end{thebibliography}

\end{document}